\documentclass[reqno,11pt]{amsart}
\usepackage[utf8]{inputenc}
\usepackage[T1]{fontenc}
\usepackage{amsfonts,amsmath,amssymb}
\usepackage[numbers,sort&compress]{natbib}
\usepackage{hyperref}
\numberwithin{equation}{section}

\newcommand{\R}{\mathbb{R}}
\newcommand{\dd}{\,\mathrm{d}}
\newtheorem{theorem}{Theorem}[section]
\newtheorem{corollary}[theorem]{Corollary}
\newtheorem{lemma}[theorem]{Lemma}
\newtheorem{proposition}[theorem]{Proposition}
\newtheorem{remark}[theorem]{Remark}
\theoremstyle{definition}
\newtheorem{definition}[theorem]{Definition}
\newtheorem{assumption}[theorem]{Assumption}

\begin{document}
\title[Rotating Euler--Riesz stars]{Relative Entropy Methods for
Rotating Stars\\ in the Compressible Euler--Riesz Equations}

\author[G.-Q. Chen]{Gui-Qiang G. Chen}
\address{Gui-Qiang G. Chen:\, Mathematical Institute, University of Oxford, Oxford OX2 6GG, UK;
School of Mathematical Sciences, Fudan University, Shanghai 200433, China;
AMSS, Chinese Academy of Sciences, Beijing 100190, China}
\email{gui-qiang.chen@maths.ox.ac.uk}

\author[D.~F. Yuan]{Difan Yuan}
\address{Difan Yuan: \,School of Mathematical Sciences,  Beijing Normal University and Laboratory of Mathematics and Complex Systems,
Ministry of Education, Beijing 100875, China}
\email{yuandf@amss.ac.cn}

\dedicatory{Dedicated to Yan Guo on the Occasion of His 60th Birthday}

\begin{abstract}
We establish finite-time relative entropy estimates 
for axisymmetric finite-energy entropy solutions of the three-dimensional attractive compressible Euler--Riesz equations near a prescribed compactly supported rotating star.
The entropy solutions are assumed to have the same total mass as the steady state and 
to satisfy the prescribed angular momentum relation. 
For general pressure laws, the compatibility with the angular momentum constraint, together with the nonnegative rotational remainder, yields a cancellation of the centrifugal terms. For axisymmetric perturbations, the gradient of the perturbed Riesz potential has no azimuthal component and therefore couples only to those momentum components controlled by the modulated kinetic energy. 
The density-weighted estimates for the Riesz potential then allow us to close the relative entropy inequality. 
Under an additional higher-integrability assumption on the density, 
the result extends to the Chandrasekhar pressure law, 
yielding a corresponding stability estimate for rotating white dwarf stars.
\end{abstract}

\date{\today}

\maketitle

\thispagestyle{empty}

\section{Introduction}
We consider the three-dimensional (3-D) attractive compressible Euler-Riesz equations (CEREs)
\begin{align}\label{1.1}
\begin{cases}
\displaystyle
\partial_t \rho+\operatorname{div}\, \mathbf{M}=0,\\
\displaystyle
\partial_t \mathbf{M}+\operatorname{div}\Big(\frac{\mathbf{M}\otimes \mathbf{M}}{\rho}\Big)+\nabla p(\rho)+\rho\nabla \mathcal{W}_{\alpha}=\boldsymbol{0},\\
\end{cases}
\end{align}
for $(t,\boldsymbol{x})\in\R_+\times\R^3$, where
$\rho:\R_+\times\R^3\to\R_+$ is the density,
$\mathbf{M}:\R_+\times\R^3\to\R^3$ is the momentum, and
$p=p(\rho)$ is a barotropic pressure law.  For the relative entropy analysis,
we restrict to $\alpha\in(0,2)$ and set
$\mathcal{W}_{\alpha}=\Phi_\alpha*\rho$, where
 \begin{equation*}
 \nabla \mathcal{W}_{\alpha}(\boldsymbol{x})=(\nabla \Phi_\alpha * \rho)(\boldsymbol{x})
 \end{equation*}
and
 \begin{equation*}
 \Phi_\alpha(\boldsymbol{x})=-\frac{|\boldsymbol{x}|^{-\alpha}}{\alpha}.
 \end{equation*}

Riesz potentials provide a natural scale of long-range interactions appearing
in plasma physics, stellar and galactic dynamics, condensed-matter models, random-matrix theory, and collective phenomena such as swarming
\cite{Campa_2014,Carrillo_2017,Carrillo_2021,Lewin_2022}.  In 3-D, $\alpha=1$
corresponds to the Newtonian (Coulomb) potential and $\mathcal{W}_1$ satisfies the Poisson equation, up to normalization.  In this case, \eqref{1.1} reduces to the compressible
Euler--Poisson equations (CEPEs).  For $\alpha\ne1$, the potential satisfies the fractional elliptic equation
$$(-\Delta)^{\frac{3-\alpha}{2}}\mathcal{W}_{\alpha}=-c_{\alpha}\rho$$ in the sense of distributions, where $c_{\alpha}>0$ is a normalization constant. 
Such equations can also be applied in astrophysical models of fractional gravity. 
For example, Jeans instability for CEREs was studied in \cite{Ourabah2024}.
The lack of a local Poisson formulation and the differing singular structures of sub-Coulomb and super-Coulomb potentials present significant analytical challenges, particularly in establishing compactness and obtaining uniform estimates for nonlocal interaction terms.

 Choi-Jeong \cite{Choi_2022} established the local classical well-posedness for the CEREs
with either attractive or repulsive interactions, and for the pressureless
system in the repulsive case; they also obtained sufficient conditions for
finite-time breakdown in isentropic and isothermal flows. 
Global classical well-posedness results are perturbative
and use different mechanisms.  For isentropic flows, global solutions
have been constructed for small initial density and velocity close to a
uniform state
\cite{Danchin_2022,Choi_2025}.  
Other results show that damping can stabilize
pressureless repulsive CEREs near a constant background
\cite{Choi_2022_2}, or exploit the dispersive structure of the
3-D repulsive CEREs for small irrotational perturbations of
a uniform equilibrium \cite{Choi_2025_2}. 
More recently, for the pressureless damped CEREs with repulsive
interactions, Chi-Shou-Xu \cite{ChiShouXu2025Critical}
established the global well-posedness for small perturbations of a
positive constant state in a critical Besov framework and obtained
decay rates.

For large initial data, most recently, Carrillo-Charles-Chen-Yuan \cite{Carrillo2025}
constructed global spherically symmetric finite-energy weak solutions of the
multidimensional CEREs for both
attractive and repulsive interactions.  
The existence of weak solutions and the nonlinear stability and instability 
of steady states were subsequently extended to general pressure laws; 
see \cite{Carrillo2026}. 
Unlike the Coulomb potential, the radial Riesz potential has no local Poisson
representation, and its singular behavior near the origin must therefore be
controlled carefully. For the polytropic case, developing uniform
higher-integrability estimates and using compensated compactness yield 
the strong
convergence of the approximate solutions and exclude the mass concentration at
the origin, even for attractive super-Coulomb potentials. 
Passing to the limit also requires  the 
strong compactness of the nonlocal interaction term.  

For the gravitational CEPEs, the global existence of finite-energy weak solutions 
for large spherically symmetric initial data was first established by 
Chen-He-Wang-Yuan \cite{Chen2021} for polytropic stars 
and subsequently extended by Chen \emph{et al.} \cite{Chen_2024} to general pressure laws.
A different class of global solutions with large data was constructed by 
Had\v{z}i\'c-Jang \cite{HadzicJang2019EP}, 
who established the existence of compactly supported physical-vacuum solutions 
in 3-D, without imposing symmetry assumptions, 
that remain close to expanding affine motions.
Concerning the existence and stability of non-rotating stellar equilibria, 
Rein \cite{Rein_2003} constructed barotropic steady states as constrained energy minimizers 
and established their conditional nonlinear stability under general perturbations. 
For general equations of state, Lin-Zeng \cite{LinZeng2022Turning} established 
a turning point principle for the linear stability of non-rotating stars and determined 
the number of unstable modes along the mass-radius curve. 
Building on the sharp linear stability criterion, Lin-Wang-Zhu \cite{LinWangZhu2025} 
subsequently proved nonlinear orbital stability for the CEPEs.
Further results concern the instability and expansion of self-gravitating stars. 
Deng-Liu-Yang-Yao \cite{Deng_2002} established the linear growth of the support for positive energy 
solutions with $\gamma>\frac43$. Jang \cite{Jang2008,Jang2014} proved the nonlinear 
instability of Lane-Emden stars for $\gamma=\frac65$ and $\gamma\in(\frac65,\frac43)$. 
More recently, Cheng \emph{et al.} \cite{Cheng_2025} constructed expanding solutions near unstable 
Lane-Emden equilibria.

Of particular relevance here are Yan Guo's contributions to both the plasma
and self-gravitating CEPEs.  
For the 3-D electron-fluid model, Guo \cite{Guo1998} constructed global smooth irrotational flows with
sufficiently small initial velocity.  
Guo-Pausader \cite{GuoPausader2011} subsequently proved
global smoothness for small-amplitude irrotational ion dynamics.  
In the gravitational setting, Guo-Had\v{z}i\'c-Jang showed that, for every $\gamma\in(1,\frac{4}{3})$, there exists an
infinite-dimensional family of spherically symmetric, compactly supported
physical vacuum solutions undergoing continued gravitational collapse; see also \cite{GuoHadzicJang2021}.  
They also rigorously
constructed the isothermal Larson-Penston self-similar collapse profile
\cite{GuoHadzicJang2021LP}.  Guo-Had\v{z}i\'c-Jang-Schrecker \cite{GuoHadzicJangSchrecker2022} constructed
real-analytic Yahil-type polytropic profiles on the full self-similar radial
domain for every $\gamma\in(1,\frac{4}{3})$.  Their recent work proves nonlinear
stability of the Larson-Penston family under radially symmetric
perturbations \cite{GuoHadzicJangSchrecker2025LPStability}.

The variational construction of rotating Euler-Poisson stars was initiated by 
Auchmuty-Beals \cite{AuchmutyBeals1971}, who established the existence of axisymmetric equilibria 
as constrained energy minimizers with prescribed mass and distribution of squared specific 
angular momentum. 
Subsequently, Li \cite{Li1991} proved the existence of uniformly rotating stellar 
equilibria with prescribed mass for sufficiently small angular velocities. 
Chanillo--Li \cite{ChanilloLi1994} established \emph{a priori} bounds on the diameters 
and the number of connected components of uniformly rotating white dwarf equilibria. 
Luo-Smoller \cite{LuoSmoller2009} later constructed energy minimizers in a broader 
class of axisymmetric configurations and, under the assumption of uniqueness up to translations 
along the rotation axis, established their conditional nonlinear orbital stability with respect 
to global weak solutions. Further existence results for general equations of state, 
including models of rotating and non-rotating white dwarfs and rotating supermassive stars, 
were obtained in \cite{Luo_2008}.
Concerning linear stability, Lin-Wang \cite{LinWang2023} derived a sharp stability criterion 
for axisymmetric perturbations in the Rayleigh-stable regime and established linear instability 
when the Rayleigh stability condition fails. 
For further developments on the existence, structure, and stability of rotating stellar equilibria, 
we refer to \cite{Charles2026,FriedmanTurkington1981,Wu2016,StraussWu2020WhiteDwarfs,JangStraussWu2023Entropy,StraussWu2017,StraussWu2019,Yuan2021Rotating} and the references therein.

In \cite{LuoSmoller2009}, Luo-Smoller employed the relative entropy method 
to establish finite-time stability around steady solutions of the CEPEs with 
$W_{\rm loc}^{1,\infty}$-regularity, 
including rotating stellar equilibria. 
For the CEREs on the torus or in smooth bounded domains, Alves \emph{et al.} \cite{Alves_2024} established 
weak-strong uniqueness through a relative entropy argument, 
under the assumption that the density of the strong solution is uniformly bounded above 
and bounded away from vacuum.
Related advances in the treatment of singular nonlocal interactions have been achieved 
through modulated energy methods. Choi-Jung \cite{ChoiJung2024Modulated} established 
modulated interaction energy estimates for singular potentials. 
For the repulsive CEREs, Alves-Grafakos-Tzavaras \cite{AlvesGrafakosTzavaras2026} developed bilinear 
fractional integral estimates, obtaining higher space-time integrability of finite-energy solutions 
together with stability estimates. 
More recently, Carrillo \emph{et al.} \cite{Carrillo2026} established relative entropy estimates 
around compactly supported, non-rotating polytropic equilibria. 
Their approach relies on finite-energy bounds and, importantly, 
does not require a strictly positive lower bound on the density of the steady state.

The proof combines the finite-energy entropy framework developed in \cite{Carrillo2025} 
with the relative entropy method and the weighted Riesz potential estimates 
established in \cite{Carrillo2026}. Motivated by the analysis in \cite{Luo_2008,LuoSmoller2009}, 
we extend the finite-time stability estimates to rotating stellar equilibria governed 
by general pressure laws. 
The main analytical difficulties and the methods used to overcome them are outlined below:

\begin{enumerate}
\item[(i)]
For the compactly supported equilibrium $\bar{\rho}$, we retain the
weighted term involving $\bar V:=V_{\bar{\rho}}$ in $d_1$ to obtain direct
control of the far-field mass. An analogous
weighted term is used for non-rotating CEPEs in
\cite{LinWangZhu2025}. Under Assumption \ref{pressure}, we allow different power-law behaviors of the pressure at low and high densities.
 With $f=\rho-\bar\rho$ and $\Gamma=\{\bar\rho>0\}$, \eqref{d1-inside} gives
\[
h(\rho)-h(\bar\rho)+\bar{V}(\rho-\bar\rho)=D_h(\rho\mid\bar\rho)\qquad\text{on }\Gamma.
\]
On the vacuum region $\Gamma^c$, the corresponding integrand reduces to $h(\rho)+\bar V\rho$. 
The positivity property \eqref{V-properties} provides control of the contribution from the vacuum region, 
while the far-field estimate \eqref{V-far} controls the mass spreading towards spatial infinity.
More precisely, Lemma \ref{d1} provides control of the integrals of $|f|^2$ and $|f|^{\gamma_2}$ 
over $\Gamma\cap\{|f|\le1\}$ and $\Gamma\cap\{|f|>1\}$, respectively, 
as well as those of $\rho^{\gamma_1}$ and $\rho^{\gamma_2}$ over $\Gamma^c\cap\{\rho\le1\}$ 
and $\Gamma^c\cap\{\rho>1\}$, respectively. It also controls the weighted vacuum term
\[
\int_{\Gamma^c}  \bar V\,\rho\,\dd\boldsymbol{x}.
\]
Combining these estimates with Lemma \ref{density-norms}, we conclude that
\[
\rho-\bar\rho\in L^s(\R^3)
\qquad\text{for any }s\in(1,\underline\gamma].
\]

\smallskip
\item[(ii)]
The rotational energy $\mathcal{K}_L(\rho)$ depends nonlinearly and nonlocally on the density through 
the cylindrical mass function $m_\rho$. Its first variation at the equilibrium density $\bar\rho$ 
is represented by $\Psi_{\bar \rho}$. 
Under the hypotheses of Lemma~\ref{d2}, the condition $L'\ge0$ yields the exact nonnegative remainder
\[
\begin{aligned}
d_2(\rho,\bar\rho)&=\int_0^1(1-\lambda)\int_0^\infty\frac{L'(m_{\rho_\lambda}(r))m_f(r)^2}{r^3}\,\dd r\,\dd\lambda
\ge0,\
&\qquad\rho_\lambda=\bar\rho+\lambda f.
\end{aligned}
\]
The angular-momentum compatibility condition further yields \eqref{swirl-K},
which relates the azimuthal kinetic energy to the rotational energy functional $\mathcal{K}_L(\rho)$.

A key observation is that combining this identity with the decompositions involving $d_1$ and $d_2$ 
leads to an exact cancellation of the two linear terms containing $\Psi_{\bar \rho}$. 
Consequently, Proposition~\ref{cancellation} yields the effective energy identity \eqref{EL-effective}, 
in which the linear term involving $\rho-\bar\rho$ carries the weight
\[
\bar G=\mu+\bar{\mathcal W}_{\alpha}.
\]
This cancellation is essential for obtaining the relative energy estimates 
in the presence of rotation and nonlocal gravitational interactions.

\smallskip
\item[(iii)]
After the rotational cancellation, the source term in Lemma \ref{relative-weak} contains
\[
F=\nabla\Phi_\alpha*(\rho-\bar\rho).
\]
The explicit kinetic part of $\mathcal{E}_L$ controls only $M^r$ and $M^z$, since the azimuthal kinetic energy has been absorbed into
$\mathcal{K}_L(\rho)$.  The $L^s$ control obtained above permits the use of Lemma \ref{Mterm}; axisymmetry then gives
\[
F\cdot\boldsymbol{e}_\theta=0, \qquad\mathbf{M}\cdot F=(M^r\boldsymbol{e}_r+M^z\boldsymbol{e}_z)\cdot F.
\]
Thus, the nonlocal Riesz interaction involves only the momentum components controlled by the relative energy,
thereby overcoming the apparent lack of control over the azimuthal component.

The remaining nonlocal contribution is controlled by Proposition \ref{weighted-Riesz}, 
which provides the requisite weighted Riesz estimates in two regimes. 
Under condition \eqref{energy-threshold}, we have
\[
\int_{\R^3}\rho|F|^2\,\dd\boldsymbol{x}
\leq C\bigl(d_1+d_1^\chi\bigr),
\qquad
\chi=\frac{7-2\alpha}{3}.
\]
Alternatively, if \eqref{high-threshold} holds with $\bar p>\frac{3}{2-\alpha}$, then
\[
\int_{\R^3}\rho|F|^2\,\dd\boldsymbol{x}
\leq
C\bigl(1+\|\rho\|_{L^{\bar p}}\bigr)^2d_1.
\]
In the latter regime, estimate \eqref{KphT} ensures that the coefficient 
on the right-hand side remains uniformly bounded on $[0,T]$. 
Together, these estimates provide the necessary control of the nonlocal interaction 
term and allow the relative energy argument to be closed.

\smallskip
\item[(iv)]
The relative total energy satisfies
\[
\mathcal{H}(\rho,\mathbf{M}\mid\bar\rho,\bar{\mathbf{M}})=\mathcal{E}_L+\mathcal{Q}_\alpha,
\qquad
\mathcal{Q}_\alpha=\frac12\int_{\R^3}(\rho-\bar\rho)\,\Phi_\alpha*(\rho-\bar\rho)\,\dd\boldsymbol{x}\le0.
\]
In Case 1, Proposition~\ref{coercivity} yields \eqref{EM}, which, combined with the total energy inequality, 
gives \eqref{energy-bound}. Consequently,
\[
\mathcal{E}_L^\chi
\leq B_0^{\chi-1}\mathcal{E}_L,
\]
allowing the nonlinear term in \eqref{source-bound} to be controlled linearly by $\mathcal{E}_L$. 

\smallskip
In Case 2, \eqref{source-bound} already provides a linear bound in $\mathcal{E}_L$.

\smallskip
Finally, the flux condition \eqref{flux-condition} justifies passing to the limit as the spatial cutoff is removed 
in the entropy inequality. Combining the preceding estimates yields \eqref{Gronwall-linear}, 
from which the desired stability estimate \eqref{main-stability} follows by Gr\"{o}nwall's lemma.
\end{enumerate}

The Chandrasekhar pressure law \eqref{WD-pressure} satisfies Assumption \ref{pressure}. 
Indeed, the asymptotic behaviors \eqref{WD-small}--\eqref{WD-large} yield
\[
\gamma_1=\frac53,\qquad
\gamma_2=\frac43,\qquad
\underline\gamma=\frac43.
\]
For the Newtonian gravitational potential $(\alpha=1)$, the high-density exponent $\gamma_2=\frac43$ 
coincides with the mass-critical exponent. 
Consequently, Case 1 of Proposition \ref{weighted-Riesz} is not applicable, since
\[
\frac43<\frac95.
\]
However, Case 2 applies because
\[
\frac43>\frac65,
\]
provided that \eqref{WD-Lp} holds for some $\bar p>3$. 
Thus, under the remaining hypotheses of Theorem \ref{main}, Corollary \ref{white-dwarf} yields 
the conditional stability estimate \eqref{WD-stability} around any prescribed rotating 
white dwarf equilibrium satisfying Definition \ref{rotating-star}.

The rest of the paper is organized as follows: 
In Section \ref{preliminaries}, we develop the necessary estimates for general
pressure laws, state the fixed-mass coercivity condition, and define finite-energy entropy solutions.
Section \ref{equilibria} presents rotating Euler-Riesz stars and angular momentum compatibility.
In Section \ref{rotational-entropy}, we prove the coercivity of $d_1$, 
the nonnegativity of $d_2$, and the centrifugal cancellation.
In Section \ref{riesz-estimates}, we prove that the gradient of the perturbed Riesz potential has 
no azimuthal component and establish the weighted Riesz estimates in the two admissible regimes. 
In Section \ref{stability}, we derive the distributional relative entropy
inequality and prove Theorem \ref{main}. 
Section \ref{relative-total-energy} presents the relative energy decomposition.
Finally, in Section \ref{white-dwarfs}, we verify the assumptions for the Chandrasekhar 
pressure law and prove Corollary \ref{white-dwarf}.

\section{Preliminaries}\label{preliminaries}
In this section, we develop the necessary estimates for general pressure laws, state the fixed-mass coercivity condition,  and define finite-energy entropy solutions.

\subsection{Internal Energy and Mechanical Entropy}
Let $p=p(\rho):[0,\infty)\to[0,\infty)$ be a barotropic pressure law. Assuming that
$\frac{p(s)}{s^2}$ is integrable near $s=0$, define the specific internal energy $e$ and the
internal energy density $h$ by
\begin{equation}\label{h}
 e(\rho) :=\int_0^\rho\frac{p(s)}{s^2}\,\dd s,\qquad h(\rho):=\rho e(\rho),
 \qquad h(0):=0.
\end{equation}
The required integrability follows from Assumption \ref{pressure} below.

For $\rho>0$, direct differentiation yields
\begin{equation}\nonumber
 \rho h'(\rho)-h(\rho)=p(\rho),\qquad h''(\rho)=\frac{p'(\rho)}{\rho}.
\end{equation}
In particular, the condition $p'(\rho)>0$ implies that $h$ is strictly convex on $(0,\infty)$.

The corresponding mechanical entropy and entropy flux are given by
\begin{equation}\nonumber
 \eta(\rho,\mathbf M)=\frac{|\mathbf M|^2}{2\rho}+h(\rho), \qquad q(\rho,\mathbf M)
 = \mathbf M \big( \frac{|\mathbf M|^2}{2\rho^2} +h'(\rho)\big).
 \end{equation}
For smooth solutions of \eqref{1.1} with positive density, the mechanical entropy balance takes the form
\begin{equation}\nonumber
 \partial_t\eta(\rho,\mathbf M)+\operatorname{div}q(\rho,\mathbf M)=-\mathbf M\cdot\nabla\mathcal W_\alpha.
\end{equation}

We impose the following structural assumptions on the pressure law, allowing different power-law behaviours at low and high densities.

\begin{assumption}[General pressure law]
\label{pressure}
The pressure function $p$ satisfies
the following conditions:
\begin{enumerate}
\item[(i)]
$p\in C([0,\infty))\cap C^1((0,\infty))$, $p(0)=0$, and $p'(\rho)>0$ for $\rho>0$.

\item[(ii)]
There exist exponents $\gamma_1\in(1,2]$ and $\gamma_2\in(1,\infty)$,
constants $0<\rho_*<1<\rho^*$, and positive constants $c_0,C_0,\widetilde c_0,c_\infty,C_\infty$ such that
\begin{align}
 &c_0\rho^{\gamma_1}\leq p(\rho)\leq C_0\rho^{\gamma_1},\qquad p'(\rho)\geq\widetilde c_0\rho^{\gamma_1-1}
 &&\mbox{when $0<\rho\leq\rho_*$}, \label{p-small}\\
 &c_\infty\rho^{\gamma_2}
 \leq p(\rho)\leq C_\infty\rho^{\gamma_2}
 &&\mbox{when $\rho\geq\rho^*$}.
 \label{p-large}
\end{align}
\end{enumerate}
\end{assumption}

Define
\[
 \underline\gamma:=\min\{\gamma_1,\gamma_2\}.
\]

\begin{lemma}[Properties of the internel energy density]\label{pressure-bounds}
Under {\rm Assumption~\ref{pressure}}, the internal energy density $h$ defined in \eqref{h} is strictly
convex on $(0,\infty)$. Moreover, there exist positive constants $c_0',C_0',c_\infty',C_\infty'$ such that
\begin{align}
 &c_0'\rho^{\gamma_1}\leq h(\rho)\leq C_0'\rho^{\gamma_1},
\quad\,\,
 0\leq h'(\rho)
 \leq C_0'\rho^{\gamma_1-1}
\qquad\quad\,\,\, \mbox{when $0\leq\rho\leq1$},
\label{pressure-small}\\
 &c_\infty'\rho^{\gamma_2}
 \leq h(\rho) \leq C_\infty'\rho^{\gamma_2},
\quad\,\,
 0\leq h'(\rho)
 \leq C_\infty'\rho^{\gamma_2-1}
\qquad\,\,\mbox{when $\rho\geq1$}.
\label{pressure-large}
\end{align}
In addition, for any $R>0$,
\begin{equation}\label{h-convexity}
 \inf_{0<s\leq R}h''(s)=\inf_{0<s\leq R}\frac{p'(s)}{s}>0.
\end{equation}
\end{lemma}

\begin{proof}
Since $p(0)=0$ and $p'(\rho)>0$ for $\rho>0$, we see that $p(\rho)>0$ for $\rho>0$. Consequently, $e$, $h$, and $h'$ are strictly positive on $(0,\infty)$.

We first establish the estimates at low densities. For $0<\rho\leq\rho_*$,
\eqref{p-small} yields
\[
 \frac{c_0}{\gamma_1-1}\rho^{\gamma_1-1} \leq e(\rho)
 = \int_0^\rho\frac{p(s)}{s^2}\,\dd s \leq\frac{C_0}{\gamma_1-1}\rho^{\gamma_1-1}.
\]
Using
$h(\rho)=\rho e(\rho)$ and $h'(\rho)=e(\rho)+\frac{p(\rho)}{\rho}$,
we obtain
\[
 c\rho^{\gamma_1} \leq h(\rho)\leq C\rho^{\gamma_1}, \quad\,\,
 0<h'(\rho)\leq C\rho^{\gamma_1-1}\qquad\,\,\mbox{for $0<\rho\leq\rho_*$},
\]
with positive  constants $c$ and $C$.

Since $\gamma_1>1$, it follows that $h'(\rho)\to0$ as $\rho\downarrow0$.
We therefore extend $h'$ continuously to the origin by setting $h'(0)=0$.
On the compact interval $[\rho_*,1]$, the functions
$\frac{h(\rho)}{\rho^{\gamma_1}}$ and $\frac{h'(\rho)}{\rho^{\gamma_1-1}}$
are continuous and strictly positive. Combining their bounds on this interval
with the preceding estimates proves \eqref{pressure-small}.

We next consider the high-density regime. For $\rho\geq\rho^*$, write
\[
 e(\rho)=e(\rho^*)+\int_{\rho^*}^{\rho}\frac{p(s)}{s^2}\,\dd s.
\]
By \eqref{p-large},
\[
\begin{aligned}
 e(\rho^*)+\frac{c_\infty}{\gamma_2-1}\bigl(\rho^{\gamma_2-1}-(\rho^*)^{\gamma_2-1}\bigr)\leq e(\rho)\leq e(\rho^*)
 +\frac{C_\infty}{\gamma_2-1}\bigl(\rho^{\gamma_2-1}-(\rho^*)^{\gamma_2-1}\bigr).
\end{aligned}
\]
Since $\gamma_2>1$, there exist $R_0\geq\rho^*$ and positive constants
$c,C$ such that
\[
 c\rho^{\gamma_2-1}\leq e(\rho)\leq C\rho^{\gamma_2-1}\qquad \mbox{when $\rho\geq R_0$}.
\]
Together with
$h(\rho)=\rho e(\rho)$ and $h'(\rho)=e(\rho)+\frac{p(\rho)}{\rho}$,
this gives
\[
 c\rho^{\gamma_2}\leq h(\rho)\leq C\rho^{\gamma_2},\quad\,\, 0<h'(\rho)\leq C\rho^{\gamma_2-1}  
\qquad\,\,\mbox{for $\rho\geq R_0$}.
\]

On the compact interval $[1,R_0]$,
the functions $\frac{h(\rho)}{\rho^{\gamma_2}}$ and $\frac{h'(\rho)}{\rho^{\gamma_2-1}}$
are continuous and strictly positive. Adjusting the constants therefore establishes \eqref{pressure-large} for all $\rho\geq1$.

Finally, the identity
\[
 h''(\rho)=\frac{p'(\rho)}{\rho}>0 \qquad \mbox{when $\rho>0$},
\]
shows that $h$ is strictly convex.

To establish \eqref{h-convexity}, observe that, for $0<s\leq\rho_*$,
\[
 h''(s)=\frac{p'(s)}s\geq\widetilde c_0s^{\gamma_1-2}\geq\widetilde c_0\rho_*^{\gamma_1-2}>0,
\]
where the last inequality uses $\gamma_1\leq2$.

If $R>\rho_*$, the continuity and strict positivity of $\frac{p'(s)}{s}$ on
$[\rho_*,R]$ imply that
\[
 \min_{\rho_*\leq s\leq R}\frac{p'(s)}s>0.
\]
Combining this bound with the low-density
estimate gives
\[
 \inf_{0<s\leq R}h''(s)>0.
\]
The case $R\leq\rho_*$ follows directly from the low-density bound. This completes the proof.
\end{proof}

\subsection{Coercivity Condition}
For a nonnegative density $\rho$ with prescribed total mass $M$, define the attractive interaction
energy by
\begin{equation}\nonumber
\mathcal P_\alpha(\rho):=-\frac12\int_{\R^3}\rho(\Phi_\alpha*\rho)\,\dd\boldsymbol{x}
=\frac{1}{2\alpha}\iint_{\R^3\times\R^3}\frac{\rho(\boldsymbol{x})\rho(\boldsymbol{y})}
{|\boldsymbol{x}-\boldsymbol{y}|^\alpha}\,\dd\boldsymbol{x}\dd\boldsymbol{y}.
\end{equation}
The corresponding total energy is
\begin{equation}\label{total-energy}
\mathcal H(\rho,\mathbf M)=\int_{\R^3}\big(\frac{|\mathbf M|^2}{2\rho}+h(\rho)\big)\,\dd\boldsymbol{x}-\mathcal P_\alpha(\rho).
\end{equation}

We assume that the solutions under consideration satisfy the total energy inequality
\begin{equation}\label{energy-ineq}
\mathcal H(\rho(t),\mathbf M(t))\leq\mathcal H(\rho_0,\mathbf M_0) \qquad\text{for {\it a.e.} }t\in(0,T).
\end{equation}

To obtain effective control of the internal and kinetic energies, we impose the following coercivity condition at the prescribed mass $M$: there exist constants
$\vartheta_M\in(0,1)$ and $C_M<\infty$ such that
\begin{equation}\label{EM}
\mathcal P_\alpha(\varrho)\leq\vartheta_M \int_{\R^3}h(\varrho)\,\dd\boldsymbol{x}
+C_M
\end{equation}
for every nonnegative density $\varrho$ satisfying
\[
\int_{\R^3}\varrho\,\dd\boldsymbol{x}=M, \qquad h(\varrho)\in L^1(\R^3).
\]

The next proposition provides sufficient conditions for \eqref{EM}, distinguishing the mass-subcritical and mass-critical regimes.

Let $C_\alpha^*$ denote the optimal constant in the interpolation form of the Hardy--Littlewood--Sobolev inequality
\begin{equation}\label{HLS-critical}
\iint_{\R^3\times\R^3}
\frac{\varrho(\boldsymbol{x}) \varrho(\boldsymbol{y})}{|\boldsymbol{x}-\boldsymbol{y}|^\alpha}
\,\dd\boldsymbol{x}\dd\boldsymbol{y}
\leq C_\alpha^*\|\varrho\|_{L^1}^{\frac{3-\alpha}{3}}\|\varrho\|_{L^{\frac{3+\alpha}{3}}}^{\frac{3+\alpha}{3}}.
\end{equation}

\begin{proposition}[Sufficient conditions for coercivity]
\label{coercivity}
Let
\[
q_c:=\frac{3+\alpha}{3}.
\]
Under {\rm Assumption \ref{pressure}}, the following statements hold{\rm:}
\begin{enumerate}
\item[{\rm(i)}]
If $\gamma_2>q_c$, then \eqref{EM} holds for every prescribed mass $M>0$.

\item[{\rm(ii)}]
Suppose that $\gamma_2=q_c$ and
\[
K_\infty:=\liminf_{\rho\to\infty}\frac{h(\rho)}{\rho^{q_c}}\in(0,\infty).
\]
Then \eqref{EM} holds whenever
\[
\frac{C_\alpha^*}{2\alpha}M^{\frac{3-\alpha}{3}}<K_\infty.
\]
\end{enumerate}
\end{proposition}

\begin{proof}
Set
\[
A_M:=\frac{C_\alpha^*}{2\alpha}M^{\frac{3-\alpha}{3}}.
\]
By the pressure estimates \eqref{pressure-small}--\eqref{pressure-large}, there exists $C>0$ such that
\[
\rho^{\gamma_2}\leq C\big(h(\rho)+\rho\big)\qquad \mbox{ for $\rho\geq0$}.
\]
Hence every admissible density $\varrho$ belongs to $L^1(\R^3)\cap L^{\gamma_2}(\R^3)$.

Since $\gamma_2\geq q_c$ in both cases, interpolation gives $\varrho\in L^{q_c}(\R^3)$.
The Hardy--Littlewood--Sobolev estimate \eqref{HLS-critical} therefore yields
\begin{equation}\label{P-HLS}
\mathcal P_\alpha(\varrho) \leq A_M\|\varrho\|_{L^{q_c}}^{q_c}.
\end{equation}

\medskip
\noindent
\textbf{Case 1. $\gamma_2>q_c$}:
Since $1<q_c<\gamma_2$, for every $\delta>0$ there exists $C_\delta>0$ such that
\[
s^{q_c}\leq\delta s^{\gamma_2}+C_\delta s\qquad \mbox{ for $s\geq0$}.
\]
Integrating this inequality and using the prescribed mass gives
\[
\|\varrho\|_{L^{q_c}}^{q_c}\leq\delta\|\varrho\|_{L^{\gamma_2}}^{\gamma_2}+C_\delta M.
\]
The pressure estimates then imply
\[
\|\varrho\|_{L^{q_c}}^{q_c}\leq C\delta\int_{\R^3}h(\varrho)\,\dd\boldsymbol{x}+C_{\delta,M}.
\]
Combining this with \eqref{P-HLS},
we obtain
\[
\mathcal P_\alpha(\varrho)\leq CA_M\delta \int_{\R^3}h(\varrho)\,\dd\boldsymbol{x}+C_{\delta,M}.
\]
Choosing $\delta>0$ sufficiently small that
\[
\vartheta_M:=CA_M\delta<1
\]
establishes \eqref{EM} for every $M>0$.

\medskip
\noindent
\textbf{Case 2. $\gamma_2=q_c$}: Assume that $A_M<K_\infty$. Choose
\[
0<\varepsilon<K_\infty-A_M.
\]
By the definition of $K_\infty$, there exists $R_\varepsilon>0$ such that
\[
h(\rho)\geq(K_\infty-\varepsilon)\rho^{q_c}\qquad \mbox{for $\rho\geq R_\varepsilon$}.
\]
Since $h\geq0$ and
\[
\rho^{q_c}\leq R_\varepsilon^{q_c-1}\rho\qquad\mbox{for $0\leq\rho\leq R_\varepsilon$},
\]
we obtain the estimate
\[
h(\rho)\geq (K_\infty-\varepsilon)\rho^{q_c}-C_\varepsilon\rho \qquad \mbox{for $\rho\geq0$}.
\]
Integrating and using the mass constraint, we conclude that
\[
\|\varrho\|_{L^{q_c}}^{q_c}\leq\frac{1}{K_\infty-\varepsilon}\Big(\int_{\R^3}h(\varrho)\,\dd\boldsymbol{x}+C_\varepsilon M\Big).
\]
Substituting this into \eqref{P-HLS} gives
\[
\begin{aligned}
\mathcal P_\alpha(\varrho)
&\leq \frac{A_M}{K_\infty-\varepsilon}\int_{\R^3}h(\varrho)\,\dd\boldsymbol{x}+C_M.
\end{aligned}
\]
Since
\[
\vartheta_M:=\frac{A_M}{K_\infty-\varepsilon}<1,
\]
the coercivity estimate \eqref{EM} follows also in the mass-critical regime.
\end{proof}

\subsection{Finite-Energy Entropy Solutions}

We now introduce the notion of finite-energy entropy solutions to \eqref{1.1}.
\begin{definition}[Finite-energy entropy solutions]
\label{entropy-solution}
A measurable pair $(\rho,\mathbf M)$ on $[0,T]\times\R^3$ is called a \emph{finite-energy entropy solution} of \eqref{1.1} with initial data $(\rho_0,\mathbf M_0)$ if
\begin{align}
&\rho\geq0,\qquad \rho\in L^\infty(0,T;L^1(\R^3)),\\ 
&h(\rho)\in L^\infty(0,T;L^1(\R^3)),\qquad
\frac{\mathbf M}{\sqrt{\rho}} \in L^\infty(0,T;L^2(\R^3)).
\end{align}
Here $\mathbf M=0$ almost everywhere on $\{\rho=0\}$, and we adopt the convention
\[
\frac{\mathbf M}{\sqrt{\rho}}:=0 \qquad\text{on }\{\rho=0\}.
\]
All quantities appearing in the weak formulations below are assumed to be locally integrable. In addition, the following conditions are satisfied:

\begin{enumerate}
\item[(i)]
\emph{Continuity equation}:
For every
$\zeta\in W_c^{1,\infty}([0,T)\times\R^3)$,
\begin{equation}\label{continuity}
\begin{aligned}
&\int_0^T\!\!\int_{\R^3}\left(\rho\,\partial_t\zeta+\mathbf M\cdot\nabla\zeta\right)\,\dd\boldsymbol{x}\dd t+\int_{\R^3}\rho_0(\boldsymbol{x})\zeta(0,\boldsymbol{x})\,\dd\boldsymbol{x}
=0.
\end{aligned}
\end{equation}

\item[(ii)]
\emph{Momentum equations}:
For every
$\boldsymbol{\psi}\in C_c^1([0,T)\times\R^3;\R^3)$,
 \begin{align}
&\int_0^T\!\!\int_{\R^3}
\Big(\mathbf M\cdot\partial_t\boldsymbol{\psi}+\frac{\mathbf M\otimes\mathbf M}{\rho}
:\nabla\boldsymbol{\psi}+p(\rho)\operatorname{div}\boldsymbol{\psi}\Big)\,\dd\boldsymbol{x}\dd t+
\int_{\R^3}\mathbf M_0\cdot\boldsymbol{\psi}(0,\boldsymbol{x})\,\dd\boldsymbol{x}
\notag\\
&\,\,=\int_0^T\!\!\int_{\R^3}\rho\nabla\mathcal W_\alpha\cdot\boldsymbol{\psi}\,\dd\boldsymbol{x}\dd t.
\label{momentum-weak}
\end{align}

\item[(iii)]
\emph{Mechanical entropy inequality}:
For every nonnegative
$\psi\in W_c^{1,\infty}([0,T)\times\R^3)$,
\begin{align}
&\int_0^T\!\!\int_{\R^3}\left(\eta(\rho,\mathbf M)\,\partial_t\psi+q(\rho,\mathbf M)\cdot\nabla\psi
\right)\,\dd\boldsymbol{x}\dd t+\int_{\R^3}\eta(\rho_0,\mathbf M_0)\psi(0,\boldsymbol{x})\,\dd\boldsymbol{x}\notag\\
&\,\,\geq \int_0^T\!\!\int_{\R^3}\psi\,\mathbf M\cdot\nabla\mathcal W_\alpha\,\dd\boldsymbol{x}\dd t.
\label{entropy-integral}
\end{align}

\item[(iv)]
\emph{Total energy inequality}:
The total energy $\mathcal H(\rho(t),\mathbf M(t))$ defined in \eqref{total-energy} is finite
for almost every $t\in(0,T)$ and satisfies
\begin{equation*}
\mathcal H(\rho(t),\mathbf M(t)) \leq\mathcal H(\rho_0,\mathbf M_0)\qquad\text{for {\it a.e.} }t\in(0,T),
\end{equation*}
as stated in \eqref{energy-ineq}.
\end{enumerate}
\end{definition}

\section{Rotating Equilibria}\label{equilibria}

Write $\boldsymbol{x}=(x_1,x_2,z)$,
\[
r=r(\boldsymbol{x})=\sqrt{x_1^2+x_2^2},\qquad\boldsymbol{e}_r=\frac1r(x_1,x_2,0),\qquad
\boldsymbol{e}_\theta=\frac1r(-x_2,x_1,0),\qquad \boldsymbol{e}_z=(0,0,1).
\]
A scalar function is axisymmetric if it depends only on $(r,z)$.  A pair $(\rho,\mathbf{M})$ is axisymmetric if $\rho=\rho(r,z)$ and, away from the rotation axis,
\[
\mathbf{M}=\mathbf{M}^r(r,z)\boldsymbol{e}_r+\mathbf{M}^\theta(r,z)\boldsymbol{e}_\theta+\mathbf{M}^z(r,z)\boldsymbol{e}_z.
\]
For an axisymmetric density $\varrho$,
define the mass inside the cylinder of radius $r$ by
\begin{equation}\label{cylinder-mass}
m_\varrho(r)=\int_{\{r(\boldsymbol{y})\le r\}}\varrho(\boldsymbol{y})\,\dd \boldsymbol{y}=2\pi\int_0^r\int_\R\varrho(s,z)s\,\dd z\dd s.
\end{equation}

Fix a total mass $M>0$.  Let $J:[0,M]\to\R$ be a prescribed angular momentum profile and set
\begin{equation}\label{L-def}
L(m):=J(m)^2.
\end{equation}
Throughout the paper, we assume
\begin{equation}\label{L-condition}
L\in C^1([0,M]),\qquad L(0)=0,\qquad L'(m)\ge0.
\end{equation}

The rotational kinetic energy is
\begin{equation}\label{K-L}
\mathcal{K}_L(\varrho)=\frac12\int_{\R^3}\varrho(\boldsymbol{x})\frac{L(m_\varrho(r(\boldsymbol{x})))}{r(\boldsymbol{x})^2}\,\dd \boldsymbol{x}.
\end{equation}

For an axisymmetric density $\bar\rho$, define the associated centrifugal potential by
\begin{equation}\label{Psi-def}
 \Psi_{\bar\rho}(r):=\int_r^\infty \frac{L(m_{\bar\rho}(s))}{s^3}\,\dd s\qquad\,\mbox{for $r>0$},
\end{equation}
whenever the integral is finite. For brevity,
we write $\bar\Psi:=\Psi_{\bar\rho}$.

\begin{definition}[Rotating Euler-Riesz stellar density]
\label{rotating-star}
A nonnegative density
\[
 \bar\rho\in L^1(\R^3)\cap L^\infty(\R^3)
\]
is called a \emph{compactly supported rotating
Euler-Riesz stellar density of mass $M$} if the following
conditions are satisfied:

\begin{enumerate}

\item[(i)]
\emph{Mass, symmetry, and finite rotational energy}:
The density $\bar\rho$ is axisymmetric,
\[
 \int_{\R^3}\bar\rho\,\dd\boldsymbol{x}=M,
 \qquad
 \mathcal K_L(\bar\rho)<\infty,
\]
and its positivity set $\Gamma:=\{\boldsymbol{x}\in\R^3:
 \bar\rho(\boldsymbol{x})>0\}$
is bounded.

\item[(ii)]
\emph{Euler-Lagrange equilibrium conditions}:
There exists a constant $\mu>0$ such that,
with
\[
 \bar{\mathcal W}_\alpha
 :=
 \Phi_\alpha*\bar\rho,
 \qquad
 \bar\Psi:=\Psi_{\bar\rho},
\]
the following conditions hold:
\begin{align}
&h'(\bar\rho)
 +\bar\Psi
 +\bar{\mathcal W}_\alpha
 +\mu
 =0
 &&\text{{\it a.e.} on }\Gamma,
 \label{EL-on}\\
 &\bar\Psi
 +\bar{\mathcal W}_\alpha
 +\mu
 \geq0
 &&\text{{\it a.e.} on }\Gamma^c.
 \label{EL-off}
\end{align}

\item[(iii)]
\emph{Stationary Euler-Riesz equations{\rm:}}
Define the azimuthal velocity field by
\begin{equation}\label{bar-u}
 \bar{\mathbf u}(\boldsymbol{x})
 :=
 \frac{
 J(m_{\bar\rho}(r(\boldsymbol{x})))
 }{r(\boldsymbol{x})}
 \boldsymbol e_\theta,
 \qquad
 \bar{\mathbf M}
 :=
 \bar\rho\,\bar{\mathbf u},
\end{equation}
where $r(\boldsymbol{x})$ denotes the
cylindrical radius. The pair
$(\bar\rho,\bar{\mathbf M})$ is required
to be a stationary weak solution of
\eqref{1.1}.
\end{enumerate}
\end{definition}

\begin{lemma}\label{rotation-weight}
Let $\bar\rho$ be the rotating star defined in {\rm Definition \ref{rotating-star}}. If
$\rho\in L^1_+(\R^3)$ is axisymmetric, has mass $M$, and $\mathcal K_L(\rho)<\infty$, then
\begin{equation}\label{rotation-L1}
\bar{\Psi}|\rho-\bar\rho|\in L^1(\R^3).
\end{equation}
Moreover,
\begin{equation}\label{rotation-bound}
\int_{\R^3}\rho\,\bar{\Psi}\,\dd\boldsymbol{x}\le \mathcal K_L(\rho)+C_{\bar\rho,L},
\end{equation}
where $C_{\bar\rho,L}<\infty$ depends only on $\bar\rho$ and $L$.
\end{lemma}

\begin{proof}
Set
\[
F(q):=\int_0^q L(s)\,\dd s.
\]
Since $L$ is nondecreasing, $F$ is convex. The supporting-line inequality at $m_{\bar\rho}(r)$ gives
\[
\begin{aligned}
m_\rho(r)L(m_{\bar\rho}(r))\leq F(m_\rho(r))-F(m_{\bar\rho}(r))+m_{\bar\rho}(r)L(m_{\bar\rho}(r)).
\end{aligned}
\]
By the definition of $\bar\Psi$ and Tonelli's theorem, it follows that
\[
\int_{\R^3}\rho\,\bar\Psi\,\dd\boldsymbol{x}\leq\int_0^\infty\frac{F(m_\rho(r))}{r^3}\,\dd r+C_{\bar\rho,L},
\]
where
\[
C_{\bar\rho,L}:=\int_0^\infty\frac{m_{\bar\rho}(r)L(m_{\bar\rho}(r))-F(m_{\bar\rho}(r))}{r^3}\,\dd r.
\]

We first verify that $C_{\bar\rho,L}$ is finite. Since $L'\geq0$,
\[
\begin{aligned}
0\leq qL(q)-F(q)=\int_0^q sL'(s)\,\dd s\leq\frac12\|L'\|_{L^\infty(0,M)}q^2\qquad \mbox{for $0\leq q\leq M$}.
\end{aligned}
\]
Because $\bar\rho$ is bounded and compactly supported, its cylindrical mass satisfies
\[
m_{\bar\rho}(r)=O(r^2)\quad\text{as }r\downarrow0,\qquad\,\,\, m_{\bar\rho}(r)=M
\quad\text{for all sufficiently large }r.
\]
Consequently, the integrand defining $C_{\bar\rho,L}$ is $O(r)$ near the origin and $O(r^{-3})$ at infinity. Then
\[
C_{\bar\rho,L}<\infty.
\]

We next relate the first term to the rotational energy. By the definition of $\mathcal K_L$,
\[
\mathcal K_L(\rho)=\frac12\int_0^\infty\frac{L(m_\rho(r))m_\rho'(r)}{r^2}\,\dd r.
\]
Since
\[
\frac{\dd}{\dd r}F(m_\rho(r))=L(m_\rho(r))m_\rho'(r)\qquad\mbox{for {\it a.e.} $r>0$},
\]
we have
\[
\mathcal K_L(\rho)=\frac12\int_0^\infty\frac{1}{r^2}\frac{\dd}{\dd r}F(m_\rho(r))\,\dd r.
\]

To justify integration by parts, we verify that the boundary terms vanish. Since $F(0)=0$ and $L$ is nonnegative and nondecreasing,
\[
0\leq\frac{F(m_\rho(r))}{r^2}=\frac1{r^2}\int_0^rL(m_\rho(s))m_\rho'(s)\,\dd s
\leq\int_0^r\frac{L(m_\rho(s))m_\rho'(s)}{s^2}\,\dd s\longrightarrow0
\]
as $r\to 0$, where the convergence follows from $\mathcal K_L(\rho)<\infty$.

At infinity, since $m_\rho(r)\leq M$,
\[
0\leq\frac{F(m_\rho(r))}{r^2}\leq\frac{F(M)}{r^2}\longrightarrow0\qquad\text{as }r\to\infty.
\]
Integration by parts therefore yields
\[
\mathcal K_L(\rho)=\int_0^\infty\frac{F(m_\rho(r))}{r^3}\,\dd r.
\]
Combining this identity with the convexity estimate gives
\[
\int_{\R^3}\rho\,\bar\Psi\,\dd\boldsymbol{x}\leq\mathcal K_L(\rho)+C_{\bar\rho,L},
\]
which proves \eqref{rotation-bound}.

Finally, applying the same estimate with $\rho=\bar\rho$ and using $\bar\Psi\geq0$, we obtain
\[
\int_{\R^3}|\rho-\bar\rho|\bar\Psi\,\dd\boldsymbol{x}\leq\int_{\R^3}(\rho+\bar\rho)\bar\Psi\,\dd\boldsymbol{x}
\leq\mathcal K_L(\rho)+\mathcal K_L(\bar\rho)+2C_{\bar\rho,L}.
\]
This establishes \eqref{rotation-L1}.
\end{proof}

For a time-dependent axisymmetric solution
$(\rho,\mathbf M)$, decompose the momentum
in cylindrical coordinates as
\begin{equation}\label{M-decomp}
 \mathbf M = M^r\boldsymbol e_r +M^\theta\boldsymbol e_\theta +M^z\boldsymbol e_z.
\end{equation}

\begin{definition}[Angular momentum compatibility]
\label{angular momentum}
An axisymmetric solution $(\rho,\mathbf M)$
is said to be \emph{compatible with $J$}
if
\begin{equation}\label{compatibility}
 M^\theta(t,\boldsymbol{x})=\rho(t,\boldsymbol{x})\frac{J(m_{\rho(t)}(r(\boldsymbol{x}))) }{r(\boldsymbol{x})} \qquad\mbox{for {\it a.e.}
$(t,\boldsymbol{x})\in(0,T)\times\R^3$},
\end{equation}
where $r(\boldsymbol{x})$ denotes the cylindrical radius.

Under this compatibility condition, the azimuthal kinetic energy is represented by the rotational energy functional:
\begin{equation}\label{swirl-K}
 \frac12\int_{\R^3} \frac{|M^\theta|^2}{\rho} \,\dd\boldsymbol{x} = \mathcal K_L(\rho)  
\qquad \mbox{for {\it a.e.} $t\in(0,T)$}.
 \end{equation}

\end{definition}

For subsequent estimates, introduce the
effective equilibrium potentials
\[
  \bar V(\boldsymbol{x}):=\bar\Psi(r(\boldsymbol{x}))+\bar{\mathcal W}_\alpha(\boldsymbol{x})
+\mu,\qquad\bar G(\boldsymbol{x}):=\bar{\mathcal W}_\alpha(\boldsymbol{x})+\mu.
\]
By the Euler--Lagrange conditions \eqref{EL-on}--\eqref{EL-off},
\begin{equation}\label{V-properties}
 \bar V=-h'(\bar\rho)\quad\text{{\it a.e.} on }\Gamma,\qquad\,\,\,\bar V\geq0
 \quad\text{{\it a.e.} on }\Gamma^c,
\end{equation}
where $\Gamma=\{\bar\rho>0\}$.

Since $\bar\rho\in L^1(\R^3)\cap L^\infty(\R^3)$ is compactly supported and $0<\alpha<2$, the equilibrium interaction potential satisfies
\begin{equation}\label{G-regularity}
 \bar{\mathcal W}_\alpha\in W^{1,\infty}(\R^3),\qquad\bar G\in W^{1,\infty}(\R^3),
 \qquad\nabla\bar G=\nabla\bar{\mathcal W}_\alpha.
\end{equation}
Indeed, the kernels $|\boldsymbol{x}|^{-\alpha}$ and $|\boldsymbol{x}|^{-\alpha-1}$ are
locally integrable in three dimensions when $\alpha<2$. Splitting the convolution
into the regions $|\boldsymbol{x}-\boldsymbol{y}|\leq1$ and $|\boldsymbol{x}-\boldsymbol{y}|>1$,
we obtain
\[
\begin{aligned}
 &|\bar{\mathcal W}_\alpha(\boldsymbol{x})|
 \leq C\|\bar\rho\|_{L^\infty}\int_{|z|\leq1}|z|^{-\alpha}\,\dd z+C\|\bar\rho\|_{L^1},\\
 &|\nabla\bar{\mathcal W}_\alpha(\boldsymbol{x})|
 \leq C\|\bar\rho\|_{L^\infty}\int_{|z|\leq1}|z|^{-\alpha-1}\,\dd z+C\|\bar\rho\|_{L^1}.
\end{aligned}
\]
Both bounds are uniform in $\boldsymbol{x}\in\R^3$, proving \eqref{G-regularity}.

Moreover, compact support and integrability of $\bar\rho$ imply
\[
 \bar{\mathcal W}_\alpha(\boldsymbol{x}) \longrightarrow0 \qquad\text{as }|\boldsymbol{x}|\to\infty.
\]
Since $\bar\Psi\geq0$ and $\mu>0$, there exists $R_0>0$ such that
\begin{equation}\label{V-far}
 \bar V(\boldsymbol{x}) \geq\frac{\mu}{2}\qquad\text{for }|\boldsymbol{x}|\geq R_0.
\end{equation}

\section{Positive Rotational Relative Entropy}
\label{rotational-entropy}

Let $\rho$ be a nonnegative axisymmetric density with the same total mass $M$ as $\bar\rho$.
Following \cite{Luo_2008,LuoSmoller2009}, we introduce the relative energy functionals
\begin{align}
d_1(\rho,\bar\rho)
&:=
\int_{\R^3}\big(h(\rho)-h(\bar\rho)+\bar V(\rho-\bar\rho)\big)\,\dd\boldsymbol{x},
\label{d1-def}\\
d_2(\rho,\bar\rho)
&:=
\mathcal K_L(\rho)-\mathcal K_L(\bar\rho)-\int_{\R^3}\bar{\Psi}(\rho-\bar\rho)\,\dd\boldsymbol{x}.
\label{d2-def}
\end{align}
Throughout this section, we assume that $\int_{\R^3}h(\rho)\,\dd\boldsymbol{x}<\infty$
and $\mathcal K_L(\rho)<\infty$. By Lemma \ref{rotation-weight},
$\bar{\Psi}|\rho-\bar\rho|\in L^1(\R^3)$, so $d_2(\rho,\bar\rho)$ is well defined and finite.

\subsection{Coercivity of the Pressure Term}

\begin{lemma}[Coercivity of $d_1$]
\label{d1}
Under {\rm Assumption \ref{pressure}}, there exists $c>0$, depending only on
$h$ and $\bar\rho$, such that
\begin{align}
d_1(\rho,\bar\rho)
&\geq c\Big(
\int_{\Gamma\cap\{|f|\leq1\}}|f|^2\,\dd\boldsymbol{x}+\int_{\Gamma\cap\{|f|>1\}}|f|^{\gamma_2}\,\dd\boldsymbol{x}
\notag\\
&\qquad\,\,\,\,
+\int_{\Gamma^c\cap\{\rho\leq1\}}\rho^{\gamma_1}\,\dd\boldsymbol{x}
+\int_{\Gamma^c\cap\{\rho>1\}}\rho^{\gamma_2}\,\dd\boldsymbol{x}
+\int_{\Gamma^c}\bar V\rho\,\dd\boldsymbol{x}\Big),
\label{d1-coercivity}
\end{align}
where $f:=\rho-\bar\rho$ and $\Gamma:=\{\bar\rho>0\}$.
\end{lemma}

\begin{proof}
For $a,b\geq0$, define the relative internal energy
\[
D_h(a\mid b):=h(a)-h(b)-h'(b)(a-b).
\]
By the Euler--Lagrange identity
\eqref{V-properties}, we have
\begin{equation}\label{d1-inside}
h(\rho)-h(\bar\rho)+\bar V(\rho-\bar\rho)=D_h(\rho\mid\bar\rho)\qquad\text{{\it a.e.} on }\Gamma.
\end{equation}

Set $B:=\|\bar\rho\|_{L^\infty}$. We first consider the region
$\Gamma\cap\{|f|\leq1\}$. By Taylor's formula and \eqref{h-convexity},
\[
\begin{aligned}
D_h(\rho\mid\bar\rho)=f^2\int_0^1(1-s)h''(\bar\rho+sf)\,\dd s
\geq\frac12\inf_{0<z\leq B+1}h''(z)\,|f|^2
\geq c|f|^2.
\end{aligned}
\]
Here the integral representation is understood by continuity when the interpolation segment reaches zero.

Next, suppose that $|f|>1$. Since $0\leq\bar\rho\leq B$, the functions $h$ and $h'$ are
bounded on $[0,B]$. Moreover, \eqref{pressure-large} implies that
$h(a)\geq ca^{\gamma_2}$ for sufficiently large $a$.
Consequently, there exist $R>B+1$ and $c>0$ such that
\[
D_h(a\mid b)\geq ca^{\gamma_2}\geq c|a-b|^{\gamma_2}
\]
whenever $a\geq R$ and $0\leq b\leq B$.

On the compact set
\[
\mathcal A_R:=\big\{(a,b)\in[0,R]\times[0,B]\,:\,|a-b|\geq1\big\},
\]
the function
\[
(a,b)\longmapsto\frac{D_h(a\mid b)}{|a-b|^{\gamma_2}}
\]
is continuous and strictly positive, by the strict convexity of $h$.
It therefore has a positive minimum on $\mathcal A_R$.
Combining the two regions gives
\[
D_h(\rho\mid\bar\rho)\geq c|f|^{\gamma_2}\qquad\text{on }\Gamma\cap\{|f|>1\}.
\]

On the vacuum region $\Gamma^c$, we have $\bar\rho=0$ and
$\bar V\geq0$ by \eqref{V-properties}. Since $h(0)=h'(0)=0$,
\[
h(\rho)-h(\bar\rho)+\bar V(\rho-\bar\rho)=h(\rho)+\bar V\rho.
\]
The pressure estimates \eqref{pressure-small}--\eqref{pressure-large} therefore yield
\[
\begin{aligned}
h(\rho)+\bar V\rho\geq c\big(\rho^{\gamma_1}\mathbf 1_{\{\rho\leq1\}}+\rho^{\gamma_2}\mathbf 1_{\{\rho>1\}}
+\bar V\rho\big).
\end{aligned}
\]
Integrating these estimates over $\Gamma$ and $\Gamma^c$ proves
\eqref{d1-coercivity}.
\end{proof}

\begin{lemma}[Density norms controlled by $d_1$]
\label{density-norms}
Let $s\in(1,\underline\gamma]$ for $\underline\gamma=\min\{\gamma_1,\gamma_2\}$ and write $d_1:=d_1(\rho,\bar\rho)$. Then
\begin{equation}\label{fp-explicit}
\|\rho-\bar\rho\|_{L^s(\R^3)}\leq
C_s\big(d_1^{\frac{1}{2}}+d_1^{\frac{1}{\gamma_1}}+d_1^{\frac{1}{s}}\big).
\end{equation}
In particular, if $0\leq d_1\leq1$, then
\begin{equation}\label{fp-small}
\|\rho-\bar\rho\|_{L^s(\R^3)}\leq C_s d_1^{\frac{1}{2}}.
\end{equation}
\end{lemma}

\begin{proof}
Set $f:=\rho-\bar\rho$. Since
\[
1<s\leq\underline\gamma \leq\gamma_1\leq2,\qquad s\leq\gamma_2,
\]
we estimate $f$ separately on the regions appearing in \eqref{d1-coercivity}.

First, since $\Gamma$ is bounded, H\"older's inequality and the quadratic part of \eqref{d1-coercivity} give
\[
\|f\|_{L^s(\Gamma\cap\{|f|\leq1\})}\leq|\Gamma|^{\frac{1}{s}-\frac{1}{2}}\|f\|_{L^2(\Gamma\cap\{|f|\leq1\})}
\leq C d_1^{\frac{1}{2}}.
\]

On $\Gamma\cap\{|f|>1\}$, the inequality $s\leq\gamma_2$ implies $|f|^s\leq|f|^{\gamma_2}$.
Then
\[
\|f\|_{L^s(\Gamma\cap\{|f|>1\})}\leq C d_1^{\frac{1}{s}}.
\]

We next consider the vacuum region $\Gamma^c$, where $f=\rho$. On $\Gamma^c\cap\{\rho>1\}$,
\[
\|\rho\|_{L^s(\Gamma^c\cap\{\rho>1\})}^s=\int_{\Gamma^c\cap\{\rho>1\}}\rho^s\,\dd\boldsymbol{x}
\leq\int_{\Gamma^c\cap\{\rho>1\}}\rho^{\gamma_2}\,\dd\boldsymbol{x}
\leq C d_1.
\]
Therefore,
\[
\|\rho\|_{L^s(\Gamma^c\cap\{\rho>1\})}\leq C d_1^{\frac{1}{s}}.
\]

\smallskip
It remains to estimate the low-density region $\Gamma^c\cap\{\rho\leq1\}$. Let $R_0$ be as in \eqref{V-far}.
On the bounded set $B_{R_0}$, H\"older's inequality and $s\leq\gamma_1$ give
\[
\|\rho\|_{L^s(\Gamma^c\cap\{\rho\leq1\}\cap B_{R_0})}
\leq |B_{R_0}|^{\frac{1}{s}-\frac{1}{\gamma_1}}
\|\rho\|_{L^{\gamma_1}(\Gamma^c\cap\{\rho\leq1\}\cap B_{R_0})}
\leq C d_1^{\frac{1}{\gamma_1}}.
\]

Finally, set
\[
A:=\Gamma^c\cap\{\rho\leq1\}\cap B_{R_0}^{\,c}.
\]
By \eqref{V-far} and \eqref{d1-coercivity},
\[
\|\rho\|_{L^1(A)}\leq C d_1,\qquad\|\rho\|_{L^{\gamma_1}(A)}\leq C d_1^{\frac{1}{\gamma_1}}.
\]
For $1<s<\gamma_1$, choose $\theta\in(0,1)$ such that
\[
\frac1s=\theta+\frac{1-\theta}{\gamma_1}.
\]
Interpolation yields
\[
\|\rho\|_{L^s(A)}\leq\|\rho\|_{L^1(A)}^\theta\|\rho\|_{L^{\gamma_1}(A)}^{1-\theta}
\leq C d_1^{\theta+\frac{1-\theta}{\gamma_1}}
=C d_1^{\frac{1}{s}}.
\]
For $s=\gamma_1$, the same estimate follows directly from \eqref{d1-coercivity}.

Combining the estimates on all regions proves \eqref{fp-explicit}.

If $0\leq d_1\leq1$, then
$d_1^{\frac{1}{\gamma_1}}\leq d_1^{\frac{1}{2}}$ and $d_1^{\frac{1}{s}}\leq d_1^{\frac{1}{2}}$,
because $s\leq\gamma_1\leq2$. Therefore, \eqref{fp-explicit} reduces to \eqref{fp-small}.
\end{proof}

\subsection{Convexity of the Rotational Energy}
We impose the following condition at the rotation axis:
\begin{equation}\label{axis-condition}
\lim_{r\downarrow0}
\frac{
\max\{L(m_\rho(r)),L(m_{\bar\rho}(r))\}
|m_{\rho-\bar\rho}(r)|
}{r^2}
=0.
\end{equation}
This condition ensures that the boundary terms arising
from the first variation of the rotational energy
vanish at the axis. The following argument is
motivated by \cite[Lemma~3.1]{LuoSmoller2009}.

\begin{lemma}[Positivity of the rotational remainder]
\label{d2}
Let $\bar\rho$ be a rotating Euler-Riesz stellar density
satisfying {\rm Definition~\ref{rotating-star}}.
Suppose that $\rho\in L^1_+(\R^3)$ is
axisymmetric, has the same mass $M$ as
$\bar\rho$, and satisfies
\[
\mathcal K_L(\rho)<\infty.
\]
Assume further that \eqref{L-condition}
and \eqref{axis-condition} hold.
Then
\begin{equation}\label{d2-positive}
d_2(\rho,\bar\rho)\geq0.
\end{equation}
More precisely, setting
\[
f:=\rho-\bar\rho,\qquad \rho_\lambda:=\bar\rho+\lambda f,\qquad m_f:=m_\rho-m_{\bar\rho},
\]
then
\begin{equation}\label{d2-variation}
\begin{aligned}
d_2(\rho,\bar\rho)
&=\int_0^1(1-\lambda)\int_0^\infty\frac{L'(m_{\rho_\lambda}(r))m_f(r)^2}{r^3}\,\dd r\,\dd\lambda.
\end{aligned}
\end{equation}
\end{lemma}

\begin{proof}
We divide the proof into five steps.

\smallskip
\noindent
\textbf{1.} {\it Cylindrical mass representation}:
For every nonnegative axisymmetric density $\varrho\in L^1(\R^3)$ of mass $M$, the cylindrical mass $m_\varrho$ is absolutely continuous and satisfies
\begin{equation}\label{mass-derivative}
m_\varrho'(r)=2\pi r\int_\R\varrho(r,z)\,\dd z\qquad\text{for {\it a.e.} }r>0.
\end{equation}
In particular,
\[
m_\varrho(0)=0,\qquad 0\leq m_\varrho(r)\leq M,\qquad \lim_{r\to\infty}m_\varrho(r)=M.
\]
Using cylindrical coordinates and \eqref{mass-derivative}, we obtain
\[
\mathcal K_L(\varrho)=\frac12\int_0^\infty\int_\R\varrho(r,z)\frac{L(m_\varrho(r))}{r^2}
\,2\pi r\,\dd z\,\dd r
=\frac12\int_0^\infty\frac{L(m_\varrho(r))m_\varrho'(r)}{r^2}\,\dd r.
\]

\smallskip
\noindent
\textbf{2.} {\it First variation on a truncated interval}:
By linearity of the cylindrical mass,
\[
m_\lambda:=m_{\rho_\lambda}=(1-\lambda)m_{\bar\rho}+\lambda m_\rho,
\]
and hence
\[
\partial_\lambda m_\lambda=m_f,\qquad \partial_\lambda m_\lambda'=m_f'.
\]

Fix $0<\varepsilon<R<\infty$ and define
\[
\mathcal K_{\varepsilon,R}(\lambda):=\frac12\int_\varepsilon^R\frac{L(m_\lambda(r))m_\lambda'(r)}{r^2}\,\dd r.
\]
Since
\[
|m_f|\leq M,\qquad |m_f'|\leq m_\rho'+m_{\bar\rho}',\qquad 0\leq m_\lambda' \leq m_\rho'+m_{\bar\rho}',
\]
and $L,L'$ are bounded on $[0,M]$, differentiation under the integral sign is justified on $[\varepsilon,R]$.

Moreover,
\[
\partial_\lambda \bigl(L(m_\lambda)m_\lambda'\bigr)
=L'(m_\lambda)m_fm_\lambda' +L(m_\lambda)m_f'
=
\frac{\dd}{\dd r}
\bigl(L(m_\lambda)m_f\bigr).
\]
Integrating first in $\lambda$ and then by parts in $r$, we obtain
\begin{align}
\mathcal K_{\varepsilon,R}(1)-\mathcal K_{\varepsilon,R}(0)
=\frac12\int_0^1\Big[\frac{L(m_\lambda(r))m_f(r)}{r^2}\Big]_{r=\varepsilon}^{r=R}\,\dd\lambda
+\int_0^1\int_\varepsilon^R\frac{L(m_\lambda(r))m_f(r)}{r^3}\,\dd r\,\dd\lambda.
\label{K-difference}
\end{align}

\smallskip
\noindent
\textbf{3.} {\it Vanishing of the boundary terms}:
Since the integrands defining the rotational energies are nonnegative, monotone convergence gives
\[
\mathcal K_{\varepsilon,R}(1)\longrightarrow\mathcal K_L(\rho),\quad\mathcal K_{\varepsilon,R}(0)
\longrightarrow\mathcal K_L(\bar\rho)\qquad\,\,\, \mbox{as $\varepsilon\to0$ and $R\to\infty$.}
\]

At infinity, the bounds $0\leq m_\lambda(R)\leq M$ and $|m_f(R)|\leq M$ imply
\[
\sup_{\lambda\in[0,1]}\Big|\frac{L(m_\lambda(R))m_f(R)}{R^2}\Big|\leq\frac{M\|L\|_{L^\infty(0,M)}}{R^2}\longrightarrow0.
\]

At the rotation axis, $m_\lambda(\varepsilon)$ lies between $m_\rho(\varepsilon)$ and
$m_{\bar\rho}(\varepsilon)$. Since $L$ is nonnegative and nondecreasing,
\[
L(m_\lambda(\varepsilon))\leq \max\{L(m_\rho(\varepsilon)),\,L(m_{\bar\rho}(\varepsilon))\}.
\]
Consequently, \eqref{axis-condition} gives
\[
\sup_{\lambda\in[0,1]}\Big|\frac{L(m_\lambda(\varepsilon))m_f(\varepsilon)}{\varepsilon^2}
\Big|
\leq\frac{\max\{L(m_\rho(\varepsilon)),L(m_{\bar\rho}(\varepsilon))\}|m_f(\varepsilon)|}{\varepsilon^2}\longrightarrow0.
\]
Thus the boundary contribution in \eqref{K-difference} vanishes in the limit.

\smallskip
\noindent
\textbf{4.} {\it Identification of the first variation}:
By \eqref{mass-derivative},
\[
|m_f'(r)|\leq 2\pi r\int_\R|f(r,z)|\,\dd z.
\]
Therefore, Lemma \ref{rotation-weight} yields
\begin{equation}\label{mf-Psi}
\int_0^\infty|m_f'(r)|\bar{\Psi}(r)\,\dd r\leq\int_{\R^3}|f(\boldsymbol{x})|\bar{\Psi}(r(\boldsymbol{x}))
\,\dd\boldsymbol{x}<\infty.
\end{equation}

Set
\[
A(s):=\frac{L(m_{\bar\rho}(s))}{s^3},\qquad \bar{\Psi}(r)=\int_r^\infty A(s)\,\dd s.
\]
Since $m_f(0)=0$ and $m_f$ is absolutely continuous,
\[
|m_f(s)|\leq\int_0^s|m_f'(r)|\,\dd r.
\]
Tonelli's theorem and \eqref{mf-Psi} therefore give
\[
\int_0^\infty A(s)|m_f(s)|\,\dd s\leq\int_0^\infty A(s)\int_0^s|m_f'(r)|\,\dd r\,\dd s
=\int_0^\infty|m_f'(r)|\int_r^\infty A(s)\,\dd s\,\dd r
<\infty.
\]

Fubini's theorem is thus applicable, and we obtain
\begin{align}
\int_{\R^3}f(\boldsymbol{x})\bar{\Psi}(r(\boldsymbol{x}))\,\dd\boldsymbol{x}
=\int_0^\infty
m_f'(r)\bar{\Psi}(r)\,\dd r
=\int_0^\infty\frac{L(m_{\bar\rho}(r))m_f(r)}{r^3}\,\dd r.
\label{K-prime-Psi}
\end{align}

\smallskip
\noindent
\textbf{5.} {\it  Convexity and the exact remainder}:
Subtracting the truncated integral
corresponding to \eqref{K-prime-Psi}
from \eqref{K-difference}, we obtain
\begin{align}
&\mathcal K_{\varepsilon,R}(1)-\mathcal K_{\varepsilon,R}(0)-\int_\varepsilon^R
\frac{L(m_{\bar\rho}(r))m_f(r)}{r^3}\,\dd r\notag\\
&\quad=\frac12\int_0^1\Big[\frac{L(m_\lambda(r))m_f(r)}{r^2}\Big]_{r=\varepsilon}^{r=R}
\,\dd\lambda
+\int_0^1\int_\varepsilon^R\frac{\big(L(m_\lambda(r))-L(m_{\bar\rho}(r))\big)
m_f(r)}{r^3}\,\dd r\,\dd\lambda.
\label{d2-truncated}
\end{align}

Since $m_\lambda-m_{\bar\rho}=\lambda m_f$ and $L$ is nondecreasing,
\begin{equation}\label{bregman-positive}
\big(L(m_\lambda)-L(m_{\bar\rho})\big)m_f\geq0.
\end{equation}
The left-hand side of \eqref{d2-truncated} converges to $d_2(\rho,\bar\rho)$ by
Steps 3--4, while the boundary term vanishes by Step 3.
The monotone convergence theorem therefore yields
\begin{equation}\label{d2-bregman}
\begin{aligned}
d_2(\rho,\bar\rho)
&=\int_0^1\int_0^\infty\frac{\big(L(m_\lambda(r))-L(m_{\bar\rho}(r))\big)m_f(r)}{r^3}\,\dd r\,\dd\lambda.
\end{aligned}
\end{equation}

Introduce the primitive
\[
\mathcal A(q):=\int_0^qL(s)\,\dd s.
\]
Then
\[
\mathcal A'(q)=L(q),\qquad \mathcal A''(q)=L'(q)\geq0,
\]
so $\mathcal A$ is convex on $[0,M]$.

For each fixed $r>0$, the fundamental theorem of calculus gives
\[
\int_0^1\big(L(m_\lambda)-L(m_{\bar\rho})\big)m_f\,\dd\lambda
=\mathcal A(m_\rho)-\mathcal A(m_{\bar\rho})-L(m_{\bar\rho})m_f.
\]
Consequently,
\[
d_2(\rho,\bar\rho)=\int_0^\infty\frac{1}{r^3}\big(\mathcal A(m_\rho(r))-\mathcal A(m_{\bar\rho}(r))
-\mathcal A'(m_{\bar\rho}(r))m_f(r)\big)\,\dd r.
\]
The integrand is nonnegative by convexity of $\mathcal A$, proving
\eqref{d2-positive}.

Finally,
\[
L(m_\lambda)-L(m_{\bar\rho})=\int_0^\lambda L'(m_\sigma)m_f\,\dd\sigma.
\]
Substituting into \eqref{d2-bregman} and applying Tonelli's theorem gives
\[
\begin{aligned}
d_2(\rho,\bar\rho)
&=
\int_0^1\int_0^\infty\int_0^\lambda\frac{L'(m_\sigma(r))m_f(r)^2}{r^3}\,\dd\sigma\,\dd r\,\dd\lambda\\
&=\int_0^1(1-\sigma)\int_0^\infty\frac{L'(m_\sigma(r))m_f(r)^2}{r^3}\,\dd r\,\dd\sigma.
\end{aligned}
\]
Relabeling $\sigma$ as $\lambda$ establishes \eqref{d2-variation}.
\end{proof}

\subsection{Modulated Entropy and Rotational Cancellation}

We now introduce the modulated entropy associated
with a rotating Euler-Riesz equilibrium.

\begin{definition}[Rotating relative entropy]
\label{EL}
For an axisymmetric state $(\rho,\mathbf M)$
satisfying the angular momentum compatibility
condition \eqref{compatibility}, define
\begin{equation}\label{EL-def}
\mathcal E_L(t)
:=
\frac12\int_{\R^3}
\Big(\frac{|M^r|^2+|M^z|^2}{\rho}\Big)
(t,\boldsymbol{x})\,\dd\boldsymbol{x}
+d_1(\rho(t),\bar\rho)
+d_2(\rho(t),\bar\rho).
\end{equation}
\end{definition}

Under the hypotheses of
Lemmas~\ref{d1} and \ref{d2},
including the axis condition
\eqref{axis-condition}, we have
\[
\mathcal E_L(t)\geq0.
\]
The following proposition establishes
an exact cancellation of the centrifugal
potential terms.

\begin{proposition}[Rotational cancellation]
\label{cancellation}
Let $(\rho,\mathbf M)$ be an axisymmetric
finite-energy state satisfying
\eqref{compatibility}. Assume that
$d_1(\rho,\bar\rho)$ and
$d_2(\rho,\bar\rho)$ are finite.
Then all terms in the following identity
are well defined and finite{\rm:}
\begin{equation}\label{EL-effective}
\begin{aligned}
\mathcal E_L(t)
&=
\int_{\R^3}
\Big(
\frac{|\mathbf M|^2}{2\rho}
+h(\rho)-h(\bar\rho)
+\bar G(\rho-\bar\rho)
\Big)
(t,\boldsymbol{x})\,\dd\boldsymbol{x}
-\mathcal K_L(\bar\rho).
\end{aligned}
\end{equation}
\end{proposition}

\begin{proof}
By angular momentum compatibility
and \eqref{swirl-K},
\[
\frac12\int_{\R^3}
\frac{|M^r|^2+|M^z|^2}{\rho}
\,\dd\boldsymbol{x}
+\mathcal K_L(\rho)
=
\frac12\int_{\R^3}
\frac{|\mathbf M|^2}{\rho}
\,\dd\boldsymbol{x}.
\]

On the other hand, adding
\eqref{d1-def} and \eqref{d2-def}
gives
\[
d_1(\rho,\bar\rho)+d_2(\rho,\bar\rho)
=
\int_{\R^3}
\big(
h(\rho)-h(\bar\rho)
+
(\bar V-\bar{\Psi})(\rho-\bar\rho)
\big)\,\dd\boldsymbol{x}
+
\mathcal K_L(\rho)
-\mathcal K_L(\bar\rho).
\]
Since
$\bar{V}-\bar{\Psi}=\bar G$,
the two linear terms involving
$\bar{\Psi}$ cancel exactly.
Combining the preceding identities
yields \eqref{EL-effective}.
\end{proof}

\section{Weighted Estimates for the Riesz Potential}\label{riesz-estimates}
We first show that the gradient of the Riesz potential generated by an axisymmetric density
perturbation has no azimuthal component. This observation is important because the
rotational relative entropy controls only the radial and axial components of momentum.

\begin{lemma}[Meridional structure of the Riesz force]
\label{Mterm}
Let $\alpha\in(0,2)$, and let $f\in L^s(\R^3)$ be axisymmetric for some
\[
1<s<\frac{3}{2-\alpha}.
\]
Then $F:=\nabla\Phi_\alpha*f$ is well defined almost everywhere and satisfies
\begin{equation}\label{no-theta}
F=F^r(r,z)\boldsymbol e_r+F^z(r,z)\boldsymbol e_z\qquad\text{a.e. in }\R^3.
\end{equation}
Consequently, for any axisymmetric momentum field $\mathbf M$,
\begin{equation}\label{MF-meridional}
\mathbf M\cdot F=(M^r\boldsymbol e_r+M^z\boldsymbol e_z)\cdot F\qquad\text{a.e.}
\end{equation}
\end{lemma}

\begin{proof}
The Hardy-Littlewood-Sobolev inequality for the gradient Riesz potential gives
\begin{equation}\label{HLS-gradient}
\|\nabla\Phi_\alpha*f\|_{L^q(\R^3)}\leq C_{\alpha,p}\|f\|_{L^p(\R^3)}
\qquad \mbox{with $\frac1q=\frac1p-\frac{2-\alpha}{3}$},
\end{equation}
for $1<p<\frac{3}{2-\alpha}$.

We first assume that $f\in C_c^\infty(\R^3)$ is axisymmetric. Since $\Phi_\alpha$ is radial, for every rotation $R$ about the $z$-axis,
\[
(\Phi_\alpha*f)(R\boldsymbol{x})=(\Phi_\alpha*f)(\boldsymbol{x}).
\]
Differentiating gives
\[
F(R\boldsymbol{x})=RF(\boldsymbol{x}).
\]
Moreover, the potential is independent of the azimuthal angle, so
\[
F\cdot\boldsymbol e_\theta=0\qquad\text{away from the rotation axis}.
\]

For general axisymmetric $f\in L^s(\R^3)$, choose a radial mollifier $\eta$ and radial cutoffs
$\chi_n\in C_c^\infty(\R^3)$ such that $0\leq\chi_n\leq1$ and $\chi_n\to1$ pointwise. Set
\[
f_n:=\chi_n(\eta_{\frac{1}{n}}*f).
\]
Then $f_n\in C_c^\infty(\R^3)$ is axisymmetric and $f_n\longrightarrow f$ in $L^s(\R^3)$.
Let
\[
F_n:=\nabla\Phi_\alpha*f_n,\qquad \frac1q=\frac1s-\frac{2-\alpha}{3}.
\]
By \eqref{HLS-gradient},
\[
\|F_n-F\|_{L^q}\leq C_{\alpha,s}\|f_n-f\|_{L^s} \longrightarrow0.
\]
Since $F_n\cdot\boldsymbol e_\theta=0$, we have
\[
\|F\cdot\boldsymbol e_\theta\|_{L^q}\leq\|F-F_n\|_{L^q}\longrightarrow0.
\]
Thus, $F$ has no azimuthal component.

Furthermore, for every fixed rotation $R$ about the $z$-axis,
\[
\|F(R\cdot)-RF\|_{L^q}\leq\|F(R\cdot)-F_n(R\cdot)\|_{L^q}
+
\|R(F_n-F)\|_{L^q}
=2\|F_n-F\|_{L^q}\longrightarrow0.
\]
Hence, $F(R\boldsymbol{x})=RF(\boldsymbol{x})$ almost everywhere. It follows that
$F^r$ and $F^z$ depend only on $(r,z)$, proving \eqref{no-theta}.
Finally, \eqref{MF-meridional} follows from \eqref{M-decomp}.
\end{proof}

We next establish the weighted Riesz potential estimates needed to control
the nonlocal source term in the relative entropy inequality. Two regimes arise:
the first requires only finite-energy control, while the second uses an additional high-integrability bound on the density and includes the Chandrasekhar pressure
law for rotating white dwarfs.

\begin{proposition}[Weighted Riesz potential estimates]
\label{weighted-Riesz}
Let $\bar\rho$ be a rotating Euler-Riesz stellar density satisfying
{\rm Definition~\ref{rotating-star}}, and let $\rho\in L^1_+(\R^3)$ be
axisymmetric with the same mass $M$. Assume that
\[
d_1:=d_1(\rho,\bar\rho)<\infty, \qquad f:=\rho-\bar\rho.
\]

Suppose that one of the following conditions holds{\rm :}

\medskip
\noindent
\emph{Case 1 (finite-energy regime):}
\begin{equation}\label{energy-threshold}
\underline\gamma\geq\frac9{7-2\alpha}.
\end{equation}

\medskip
\noindent
\emph{Case 2 (high-integrability regime):}
\begin{equation}\label{high-threshold}
\underline\gamma\geq\max\big\{1,\frac6{7-2\alpha}\big\},\qquad\rho\in L^{\bar p}(\R^3)
\quad\text{for some }\bar p>\frac3{2-\alpha}.
\end{equation}
Then, in {\rm Case 1},
\begin{equation}\label{weightedcase1}
\int_{\R^3}\rho\,|\nabla\Phi_\alpha*f|^2\,\dd\boldsymbol{x}\leq C(d_1+d_1^\chi),\qquad \chi:=\frac{7-2\alpha}{3}\geq1.
\end{equation}

In {\rm Case 2},
\begin{equation}\label{weightedhigh}
\int_{\R^3}\rho\,|\nabla\Phi_\alpha*f|^2\,\dd\boldsymbol{x}\leq C\left(1+\|\rho\|_{L^{\bar p}}\right)^2d_1.
\end{equation}
\end{proposition}

\begin{proof}
We may assume $d_1>0$, since $d_1=0$ implies $f=0$ almost everywhere by Lemma \ref{d1}.

\medskip
\noindent
{Case 1.}
Set
\[
p_*:=\frac9{7-2\alpha},\qquad q_*:=\frac9{1+\alpha}.
\]
Then
\[
\frac1{q_*}=\frac1{p_*}-\frac{2-\alpha}{3},\qquad\frac1{p_*}+\frac2{q_*}=1.
\]
By \eqref{energy-threshold}, $p_*\leq\underline\gamma$. Choose $p_b$ such that
$\max\big\{1,\frac6{7-2\alpha}\big\}<p_b<p_*$,
and define
\[
\frac1{q_b}:=\frac1{p_b}-\frac{2-\alpha}{3}.
\]
Then
$1<p_b<p_*\leq\underline\gamma$ and $q_b>2$.

Since $\bar\rho$ is bounded and supported in the bounded set $\Gamma$, H\"older's inequality and
\eqref{HLS-gradient} yield
\[
\int_{\R^3}\bar\rho\,|\nabla\Phi_\alpha*f|^2\,\dd\boldsymbol{x}\leq\|\bar\rho\|_{L^\infty}
|\Gamma|^{1-2/q_b}\|\nabla\Phi_\alpha*f\|_{L^{q_b}}^2
\leq C\|f\|_{L^{p_b}}^2.
\]

Similarly, using $1/p_*+2/q_*=1$, we obtain
\[
\int_{\R^3}|f|\,|\nabla\Phi_\alpha*f|^2\,\dd\boldsymbol{x}
\leq\|f\|_{L^{p_*}}\|\nabla\Phi_\alpha*f\|_{L^{q_*}}^2
\leq C\|f\|_{L^{p_*}}^3.
\]
Since $\rho\leq\bar\rho+|f|$,
\[
\int_{\R^3}\rho\,|\nabla\Phi_\alpha*f|^2\,\dd\boldsymbol{x}\leq C\left(\|f\|_{L^{p_b}}^2+\|f\|_{L^{p_*}}^3\right).
\]

Applying \eqref{fp-explicit} with $s=p_b,p_*$ gives
\[
\|f\|_{L^{p_b}}^2\leq C\big(d_1+d_1^{\frac{2}{\gamma_1}}
+d_1^{\frac{2}{p_b}}\big),\qquad \|f\|_{L^{p_*}}^3
\leq C\big(d_1^{\frac{3}{2}}+d_1^{\frac{3}{\gamma_1}}+d_1^{\frac{3}{p_*}}\big).
\]
Choose $p_b$ sufficiently close to $p_*$ so  that $p_b>\frac{2p_*}{3}$.
Since $p_*\leq\underline\gamma\leq\gamma_1\leq2$, all exponents of $d_1$ above belong to $[1,\chi]$, 
where
$\chi=\frac3{p_*}=\frac{7-2\alpha}{3}$.
Using
\[
d_1^\beta\leq d_1+d_1^\chi \qquad \mbox{for $1\leq\beta\leq\chi$},
\]
we obtain \eqref{weightedcase1}.

\medskip
\noindent
{Case 2.}
Since $\bar\rho$ is bounded and compactly supported,
\[
f\in L^1(\R^3)\cap L^{\bar p}(\R^3),\qquad\|f\|_{L^1}\leq2M.
\]
Let $\bar p'$ denote the conjugate exponent of $\bar p$.
The condition $\bar p>\frac{3}{2-\alpha}$ implies
\[
(\alpha+1)\bar p'<3.
\]

Splitting the convolution into $|\boldsymbol{x}-\boldsymbol{y}|<1$ and its complement, H\"older's
inequality gives
\[
\|\nabla\Phi_\alpha*f\|_{L^\infty}
\leq C_{\alpha,\bar p}\left(\|f\|_{L^{\bar p}}+\|f\|_{L^1}\right)
\leq B_{\bar p},
\]
where
$B_{\bar p}:=C_{\alpha,\bar p}\left(2M+\|\bar\rho\|_{L^{\bar p}}+\|\rho\|_{L^{\bar p}}\right)$.

Set
\[
\ell_\alpha:=\max\big\{1,\frac6{7-2\alpha}\big\}.
\]
By \eqref{high-threshold}, we may choose
\[
\ell_\alpha\leq p_e \leq\underline\gamma,\qquad 1<p_e<\frac3{2-\alpha},
\]
sufficiently close to $\ell_\alpha$ so that, with
\[
\frac1{q_0}:=\frac1{p_e}-\frac{2-\alpha}{3},
\]
we have
\[
2\leq q_0<4.
\]
The Hardy--Littlewood--Sobolev inequality then gives
\[
\|\nabla\Phi_\alpha*f\|_{L^{q_0}}\leq C\|f\|_{L^{p_e}}.
\]

Suppose first that $0<d_1\leq1$. By \eqref{fp-small},
\[
\|\nabla\Phi_\alpha*f\|_{L^{q_0}} \leq Cd_1^{\frac{1}{2}}.
\]
Decompose
\[
|f|=g_2+g_0+g_\infty,
\]
where
$g_2:=|f|\mathbf1_{\Gamma\cap\{|f|\leq1\}}$, $g_0 :=\rho\mathbf1_{\Gamma^c\cap\{\rho\leq1\}}$,
and $g_\infty:=|f|\mathbf1_{\Gamma\cap\{|f|>1\}}+\rho\mathbf1_{\Gamma^c\cap\{\rho>1\}}$.
Set
\[
s_2:=2,\qquad s_0:=\gamma_1,\qquad s_\infty:=\min\{\gamma_2,2\}.
\]
Since $s_i\in(1,2]$, the coercivity estimate \eqref{d1-coercivity} gives
\[
\|g_i\|_{L^{s_i}}^{s_i}\leq Cd_1\qquad \mbox{for $i\in\{2,0,\infty\}$}.
\]

For each $i$, define
\[
m_i:=\frac{2s_i}{s_i-1},\qquad \theta_i:=\frac{q_0}{m_i}.
\]
Since $m_i\geq4>q_0$, interpolation between $L^{q_0}$ and $L^\infty$ yields
\[
\|\nabla\Phi_\alpha*f\|_{L^{m_i}}\leq\|\nabla\Phi_\alpha*f\|_{L^{q_0}}^{\theta_i}
\|\nabla\Phi_\alpha*f\|_{L^\infty}^{1-\theta_i}
\leq CB_{\bar p}^{1-\theta_i} d_1^{\frac{\theta_i}{2}}.
\]
By H\"older's inequality,
\[
\int_{\R^3}g_i\,|\nabla\Phi_\alpha*f|^2\,\dd\boldsymbol{x}
\leq \|g_i\|_{L^{s_i}} \|\nabla\Phi_\alpha*f\|_{L^{m_i}}^2
\leq CB_{\bar p}^{2(1-\theta_i)} d_1^{\frac{1}{s_i}+\theta_i}.
\]

Moreover,
\[
\frac1{s_i}+\theta_i-1=\frac{s_i-1}{s_i}\big(\frac{q_0}{2}-1\big)\geq0.
\]
Since $d_1\leq1$,
\[
\int_{\R^3}g_i\,|\nabla\Phi_\alpha*f|^2\,\dd\boldsymbol{x}\leq C\big(1+B_{\bar p}^2\big)d_1.
\]

For the equilibrium density,
\[
\int_{\R^3}\bar\rho\,|\nabla\Phi_\alpha*f|^2\,\dd\boldsymbol{x}\leq\|\bar\rho\|_{L^\infty}
|\Gamma|^{1-\frac{2}{q_0}}\|\nabla\Phi_\alpha*f\|_{L^{q_0}}^2\leq Cd_1.
\]
Since $\rho\leq\bar\rho+|f|$, summing these estimates gives
\[
\int_{\R^3}\rho\,|\nabla\Phi_\alpha*f|^2\,\dd\boldsymbol{x}
\leq C(1+B_{\bar p}^2)d_1\qquad \mbox{ when $0<d_1\leq1$}.
\]

Finally, if $d_1>1$, then
\[
\int_{\R^3}\rho\,|\nabla\Phi_\alpha*f|^2\,\dd\boldsymbol{x}
\leq M\|\nabla\Phi_\alpha*f\|_{L^\infty}^2
\leq MB_{\bar p}^2d_1.
\]
Combining the two ranges of $d_1$ and the definition of $B_{\bar p}$ proves \eqref{weightedhigh}.
\end{proof}

\section{Relative Entropy Inequality and Finite-Time Stability}\label{stability}
In this section,  we derive the distributional relative entropy inequality and prove Theorem \ref{main}. 

\subsection{The Effective Entropy Inequality}

We first derive the relative entropy inequality associated with the effective equilibrium potential $\bar G$.
Define the effective entropy density and flux by
\begin{align}
\eta_{\ast}(\rho,\mathbf M):=\frac{|\mathbf M|^2}{2\rho}+h(\rho)-h(\bar\rho)+\bar G\,(\rho-\bar\rho),\qquad\,\,
q_{\ast}(\rho,\mathbf M):=q(\rho,\mathbf M)+\bar G\,\mathbf M.
\notag
\end{align}
The rotational cancellation established in Proposition~\ref{cancellation} ensures that the centrifugal potential does not appear explicitly in these expressions. The remaining equilibrium potential is $\bar G=\mu+\bar{\mathcal W}_\alpha$.

\begin{lemma}[Relative entropy inequality]
\label{relative-weak}
Let $(\rho,\mathbf M)$ be an axisymmetric finite-energy entropy weak solution of \eqref{1.1} in the sense of {\rm Definition~\ref{entropy-solution}}.
Assume that
\[
\mathbf M\cdot\bigl(\nabla\Phi_\alpha*(\rho-\bar\rho)\bigr)\in L^1_{\rm loc}((0,T)\times\R^3).
\]
Then, for every nonnegative $\psi\in W_c^{1,\infty}([0,T)\times\R^3)$,
\begin{align}
&\int_0^T\!\!\int_{\R^3}\big(\eta_{\ast}(\rho,\mathbf M)\,\partial_t\psi+q_{\ast}(\rho,\mathbf M)\cdot\nabla\psi\big)
\,\dd\boldsymbol{x}\dd t
+\int_{\R^3}\eta_{\ast}(\rho_0,\mathbf M_0)\psi(0,\boldsymbol{x})\,\dd\boldsymbol{x}\notag\\
&\quad\geq\int_0^T\!\!\int_{\R^3}\psi\,\mathbf M\cdot\bigl(\nabla\Phi_\alpha*(\rho-\bar\rho)\bigr)\,\dd\boldsymbol{x}\dd t.
\label{relative-test}
\end{align}
Equivalently,
\begin{equation}\label{relative-ineq}
\partial_t\eta_{\ast}(\rho,\mathbf M)+\operatorname{div}q_{\ast}(\rho,\mathbf M)\leq-\mathbf M\cdot
\bigl(\nabla\Phi_\alpha*(\rho-\bar\rho)\bigr)
\qquad
\mbox{in $\mathcal D'((0,T)\times\R^3)$}.
\end{equation}
If $(\rho,\mathbf M)$ also satisfies the angular momentum compatibility condition \eqref{compatibility},
then, whenever the terms are finite,
\begin{equation}\label{int-effective}
\int_{\R^3}\eta_{\ast}(\rho,\mathbf M)(t,\boldsymbol{x})\,\dd\boldsymbol{x}=\mathcal E_L(t)+\mathcal K_L(\bar\rho).
\end{equation}
\end{lemma}

\begin{proof}
By \eqref{G-regularity}, $\bar G\in W^{1,\infty}(\R^3)$. Thus $\bar G\psi$ is an admissible test function
in the weak continuity equation \eqref{continuity}. Since $\bar G$ is independent of time, we obtain
\begin{align}
&\int_0^T\!\!\int_{\R^3}
\left(\bar G\rho\,\partial_t\psi+\bar G\,\mathbf M\cdot\nabla\psi\right)\,\dd\boldsymbol{x}\dd t
+\int_{\R^3}\bar G(\boldsymbol{x})\rho_0(\boldsymbol{x})\psi(0,\boldsymbol{x})\,\dd\boldsymbol{x}
\notag\\
&\quad=-\int_0^T\!\!\int_{\R^3}\psi\,\mathbf M\cdot\nabla\bar G\,\dd\boldsymbol{x}\dd t.
\label{weighted-mass}
\end{align}

Adding \eqref{weighted-mass} to the mechanical entropy inequality \eqref{entropy-integral},
and using $\nabla\bar G=\nabla\bar{\mathcal W}_\alpha$, we find that the source term becomes
\[
\mathbf M\cdot
\bigl(\nabla\mathcal W_\alpha-\nabla\bar{\mathcal W}_\alpha\bigr)
=\mathbf M\cdot\bigl(\nabla\Phi_\alpha*(\rho-\bar\rho)\bigr).
\]

Next, introduce the time-independent function
\[
C(\boldsymbol{x}):=-h(\bar\rho(\boldsymbol{x}))-\bar G(\boldsymbol{x})\bar\rho(\boldsymbol{x}).
\]
Since $\bar\rho$ is bounded and compactly supported and $\bar G$ is bounded, $C\in L^1(\R^3)$. Consequently,
\[
\int_0^T\!\!\int_{\R^3} C(\boldsymbol{x})\partial_t\psi\,\dd\boldsymbol{x}\dd t
+\int_{\R^3}C(\boldsymbol{x})\psi(0,\boldsymbol{x})\,\dd\boldsymbol{x}=0.
\]
Adding this identity gives \eqref{relative-test}. The distributional inequality
\eqref{relative-ineq} follows immediately.

Finally, under angular momentum compatibility, \eqref{int-effective} follows from the
rotational cancellation identity \eqref{EL-effective}.
\end{proof}

To pass from the local entropy inequality to a global-in-space estimate, we impose the flux integrability condition
\begin{equation}\label{flux-condition}
\rho|\mathbf u|^3+\rho\bigl(h'(\rho)\bigr)^{\frac{3}{2}}\in L^1((0,T)\times\R^3)\qquad\mbox{for $\mathbf u:=\frac{\mathbf M}{\rho}$}.
\end{equation}
Here and throughout, $\mathbf u=0$ on $\{\rho=0\}$. This condition will ensure that the entropy flux vanishes in the limit when the spatial cutoff is removed.

\subsection{Main Theorem}

We now establish finite-time stability of compactly supported rotating Euler--Riesz equilibria within the class of axisymmetric finite-energy entropy weak solutions satisfying angular momentum compatibility.

\begin{theorem}[Weighted $L^2$-stability around rotating Euler--Riesz stars]
\label{main}
Let $\alpha\in(0,2)$, let $p$ satisfy {\rm Assumption~\ref{pressure}}, and let $(\bar\rho,\bar{\mathbf M})$ represent a compactly
supported rotating Euler--Riesz star of mass $M$ in the sense of {\rm Definition~\ref{rotating-star}}.
Let $J$ be its prescribed angular momentum function, with $L=J^2$ satisfying
\eqref{L-condition}.
Let $(\rho,\mathbf M)$ be an axisymmetric finite-energy entropy solution of \eqref{1.1} on $[0,T]$, in the sense of
{\rm Definition \ref{entropy-solution}}. Assume further that the following conditions hold{\rm :}
\begin{enumerate}
\item[(i)]
The solution satisfies the angular momentum compatibility condition \eqref{compatibility}
for almost every $t\in(0,T)$. The initial data $(\rho_0,\mathbf M_0)$ satisfy the same compatibility condition and
\[
\int_{\R^3}\rho_0(\boldsymbol{x})\,\dd\boldsymbol{x}=M.
\]
Moreover, the axis condition \eqref{axis-condition} holds for a.e. $t\in(0,T)$, as well as for the initial density
$\rho_0$.

\item[(ii)]
Either the finite-energy threshold \eqref{energy-threshold} holds, or the high-integrability condition
\eqref{high-threshold} holds with
\begin{equation}\label{KphT}
K_{\bar p,T}:=\operatorname*{ess\,sup}_{0<t<T}\|\rho(t)\|_{L^{\bar p}(\R^3)}<\infty.
\end{equation}

\item[(iii)]
The flux integrability condition \eqref{flux-condition} holds.
\end{enumerate}
Then there exists a constant
\begin{equation}\label{CT-cases}
C_T=
\begin{cases}
C(T,\bar\rho,\mathcal E_L(0))
&\,\,\,\text{in {\rm Case 1}},\\
C(T,\bar\rho,K_{\bar p,T})
&\,\,\,\text{in {\rm Case 2}},
\end{cases}
\end{equation}
depending also on the fixed parameters of the problem, such that
\begin{equation}\label{main-stability}
\mathcal E_L(t)\leq C_T\mathcal E_L(0)\qquad\text{for a.e. }t\in(0,T).
\end{equation}

Moreover, define the weighted distance functional
\begin{align}
\mathcal D(t):={}& d_2(\rho(t),\bar\rho)+\frac12\int_{\R^3}\frac{|M^r|^2+|M^z|^2}{\rho}\dd\boldsymbol{x}
+\int_{\Gamma\cap\{|\rho-\bar\rho|\leq1\}}|\rho-\bar\rho|^2\,\dd\boldsymbol{x}\notag\\
&+\int_{\Gamma\cap\{|\rho-\bar\rho|>1\}}
|\rho-\bar\rho|^{\gamma_2}\,\dd\boldsymbol{x}+
\int_{\Gamma^c\cap\{\rho\leq1\}}\rho^{\gamma_1}\,\dd\boldsymbol{x}+\int_{\Gamma^c\cap\{\rho>1\}}\rho^{\gamma_2}\,\dd\boldsymbol{x}+\int_{\Gamma^c}\bar{V}\rho\,\dd\boldsymbol{x}.
\label{distance-def}
\end{align}
Then
\begin{equation}\label{distance-bound}
\mathcal D(t)\leq C\mathcal E_L(t)\leq CC_T\mathcal E_L(0)\qquad\text{for a.e. }t\in(0,T)
\end{equation}
for some constant $C>0$ independent of $T>0$.
\end{theorem}

\begin{proof}
We divide the proof into six steps.

\medskip
\noindent
\textbf{1.} {\it Finiteness and positivity of the modulated entropy}:

We first verify that $\mathcal E_L(t)$ is well defined and nonnegative for almost every $t\in(0,T)$.

By angular momentum compatibility \eqref{compatibility} and the finite kinetic energy assumption, the identity
\eqref{swirl-K} gives
\[
\mathcal K_L(\rho(t))=\frac12\int_{\R^3}\frac{|M^\theta|^2}{\rho}\,\dd\boldsymbol{x}<\infty
\qquad\mbox{for {\it a.e.} $t$}.
\]

The finite internal energy, conservation of mass, and boundedness of $\bar G$ ensure integrability of the terms involving $h(\rho)$ and $\bar G(\rho-\bar\rho)$. Moreover, Lemma~\ref{rotation-weight}, together with the finiteness of
$\mathcal K_L(\rho)$ and $\mathcal K_L(\bar\rho)$, gives
\[
\int_{\R^3}\bar{\Psi}\,|\rho-\bar\rho|\,\dd\boldsymbol{x}<\infty.
\]
Thus, $d_1(\rho,\bar\rho)$ and $d_2(\rho,\bar\rho)$ are finite.

By Lemma~\ref{d1},
\[
d_1(\rho,\bar\rho)\geq0,
\]
while Lemma~\ref{d2}, under the axis condition \eqref{axis-condition}, yields
\[
d_2(\rho,\bar\rho)\geq0.
\]
Consequently, the definition \eqref{EL-def} gives
\begin{equation}\label{EL-positive}
\mathcal E_L(t)\geq0 \qquad\text{for {\it a.e.} }t\in(0,T).
\end{equation}

\medskip
\noindent
\textbf{2}. {\it A priori bounds for the modulated entropy}:
Set
\begin{equation}\label{A-def}
A(t):=\int_{\R^3}\big(\frac{|\mathbf M|^2}{2\rho}+h(\rho)\big)(t,\boldsymbol{x})\,\dd\boldsymbol{x}.
\end{equation}
By Definition \ref{entropy-solution},
\[
A\in L^\infty(0,T).
\]

We first note that the continuity equation implies conservation of total mass. Indeed, testing \eqref{continuity} with
spatial cutoffs approaching one and using
\[
\int_{\R^3}|\mathbf M|\,\dd\boldsymbol{x}\leq\Big(\int_{\R^3}\rho\,\dd\boldsymbol{x}\Big)^{1/2}
\Big(\int_{\R^3}\frac{|\mathbf M|^2}{\rho}\,\dd\boldsymbol{x}\Big)^{1/2},
\]
we may pass to the limit in the weak continuity equation. Then
\[
\int_{\R^3}\rho(t,\boldsymbol{x})\,\dd\boldsymbol{x}
=M\qquad\text{for {\it a.e.} }t\in(0,T).
\]

By the rotational cancellation identity \eqref{EL-effective},
\[
\mathcal E_L(t)=A(t)-\int_{\R^3}h(\bar\rho)\,\dd\boldsymbol{x}+\int_{\R^3}\bar G(\rho-\bar\rho)\,\dd\boldsymbol{x}
-\mathcal K_L(\bar\rho).
\]
Since $\bar G\in L^\infty(\R^3)$ and both densities have mass $M$,
\[
\Big|\int_{\R^3}\bar G(\rho-\bar\rho)\,\dd\boldsymbol{x}\Big|\leq2M\|\bar G\|_{L^\infty}.
\]
Together with \eqref{EL-positive}, this gives
\begin{equation}\label{EL-bound}
0\leq\mathcal E_L(t)\leq A(t)+\int_{\R^3}h(\bar\rho)\,\dd\boldsymbol{x}+\mathcal K_L(\bar\rho)+2M\|\bar G\|_{L^\infty}.
\end{equation}
In particular,
\[
\mathcal E_L\in L^\infty(0,T).
\]

In Case 1, the pressure threshold implies
\[
\gamma_2\geq\underline\gamma\geq\frac9{7-2\alpha}>\frac{3+\alpha}{3}.
\]
Thus Proposition \ref{coercivity} provides the fixed-mass coercivity estimate
\eqref{EM}. Consequently,
\begin{align}
\mathcal H(\rho(t),\mathbf M(t))
&=A(t)-\mathcal P_\alpha(\rho(t))\notag\\
&\geq\int_{\R^3}\frac{|\mathbf M|^2}{2\rho}\,\dd\boldsymbol{x}+(1-\vartheta_M)\int_{\R^3}h(\rho)
\,\dd\boldsymbol{x}-C_M\notag\\
&\geq c_M A(t)-C_M,
\label{energy-lower}
\end{align}
where
$c_M:=\min\{1,1-\vartheta_M\}>0$.

Since $\mathcal P_\alpha(\rho_0)\geq0$, we have
\[
\mathcal H(\rho_0,\mathbf M_0)\leq A(0).
\]
The total energy inequality \eqref{energy-ineq} therefore yields
\begin{equation}\label{A-bound}
A(t)\leq c_M^{-1}\bigl(A(0)+C_M\bigr) \qquad\text{for {\it a.e.} }t\in(0,T).
\end{equation}

Applying \eqref{EL-effective} to the initial data, we obtain
\begin{align}
A(0)={}&\mathcal E_L(0)+\int_{\R^3}h(\bar\rho)\,\dd\boldsymbol{x}
-\int_{\R^3}\bar G(\rho_0-\bar\rho)\,\dd\boldsymbol{x}+\mathcal K_L(\bar\rho)\notag\\
\leq{}&\mathcal E_L(0)+\int_{\R^3}h(\bar\rho)\,\dd\boldsymbol{x}
+\mathcal K_L(\bar\rho)+2M\|\bar G\|_{L^\infty}.
\label{A0-bound}
\end{align}

Combining \eqref{EL-bound}, \eqref{A-bound}, and \eqref{A0-bound}, we conclude that
\begin{equation}\label{energy-bound}
\operatorname*{ess\,sup}_{0<t<T}\mathcal E_L(t)\leq B_0:=C_0\bigl(1+\mathcal E_L(0)\bigr),
\end{equation}
where $C_0$ depends only on the fixed parameters, the prescribed mass, and the equilibrium.

In Case 2, the fixed-mass coercivity estimate \eqref{EM} is not required.
The uniform boundedness of $\mathcal E_L$ already follows from the finite-energy assumption and \eqref{EL-bound}.

\medskip
\noindent
\textbf{3}. {\it Integrability of the nonlocal source term}: Set
\[
f:=\rho-\bar\rho,\qquad F:=\nabla\Phi_\alpha*f.
\]
By Lemma~\ref{Mterm}, axisymmetry implies that $F$ has no azimuthal component. Therefore,
\begin{equation}\label{sidentity}
\mathbf M\cdot F=(M^r\boldsymbol e_r+M^z\boldsymbol e_z)\cdot F\qquad\text{{\it a.e.}}
\end{equation}

Using Cauchy--Schwarz and Young's inequality,
we obtain
\begin{align}
\int_{\R^3}|\mathbf M\cdot F|\,\dd\boldsymbol{x}\leq\frac12\int_{\R^3}\frac{|M^r|^2+|M^z|^2}{\rho}\,\dd\boldsymbol{x}
+\frac12\int_{\R^3}\rho|F|^2\,\dd\boldsymbol{x}.\notag
\end{align}
Since $d_1\leq\mathcal E_L$, Proposition \ref{weighted-Riesz} yields
\begin{equation}\label{source-bound}
\int_{\R^3}|\mathbf M\cdot F|\,\dd\boldsymbol{x}\leq
\begin{cases}
C\bigl(\mathcal E_L(t)+\mathcal E_L(t)^\chi\bigr)
&\quad\text{in Case 1},\\[1mm]
C_K\mathcal E_L(t)
&\quad\text{in Case 2},
\end{cases}
\end{equation}
where
$\chi=\frac{7-2\alpha}{3}\geq1$ and $C_K=C(1+K_{\bar p,T})^2$.

In Case 1, the right-hand side is uniformly bounded by \eqref{energy-bound}.
In Case 2, the same conclusion follows from \eqref{EL-bound} and
\eqref{KphT}. Hence
\begin{equation}\label{source-L1}
\mathbf M\cdot F\in L^1((0,T)\times\R^3).
\end{equation}
In particular, the integrability hypothesis of Lemma \ref{relative-weak} is satisfied.

\medskip
\noindent
\textbf{4}. {\it Removal of the spatial cutoff}:
Choose $\zeta\in C_c^\infty(\R^3)$ such that
\[
0\leq\zeta\leq1,\qquad \zeta=1\ \text{on }B_1,\qquad \zeta=0\ \text{on }B_2^c,
\]
and define
\[
\zeta_R(\boldsymbol{x}):=\zeta(\frac{\boldsymbol{x}}{R}).
\]
Fix an increasing sequence $R_j\to\infty$.

Since $\bar G$ is bounded,
\[
|\eta_{\ast}(\rho,\mathbf M)|\leq{}\frac{|\mathbf M|^2}{2\rho}+h(\rho)+C\rho+h(\bar\rho)+C\bar\rho.
\]
The finite-energy bounds and mass conservation therefore imply $\eta_{\ast}\in L^\infty(0,T;L^1(\R^3))$.
Consequently, for each $j$,
\[
I_j(\tau):=\int_{\R^3}\zeta_{R_j}(\boldsymbol{x})\,\eta_{\ast}(\rho,\mathbf M)(\tau,\boldsymbol{x})\,\dd\boldsymbol{x}
\]
belongs to $L^1(0,T)$.

We fix $t\in(0,T)$ that is a Lebesgue point of every $I_j$. Such times form
a set of full measure, since the sequence is countable. For $0<\varepsilon<T-t$, define
\[
\vartheta_\varepsilon(\tau):=
\begin{cases}
1 & \mbox{for $0\leq\tau\leq t$},\\
1-\dfrac{\tau-t}{\varepsilon}
& \mbox{for $t<\tau<t+\varepsilon$},\\
0& \mbox{for $t+\varepsilon\leq\tau<T$}.
\end{cases}
\]
For each fixed $j$, the nonnegative
function
\[
\psi(\tau,\boldsymbol{x}):=\vartheta_\varepsilon(\tau)\zeta_{R_j}(\boldsymbol{x})
\]
is admissible in \eqref{relative-test}. Since 
$\partial_\tau\vartheta_\varepsilon=-\varepsilon^{-1}$ on $(t,t+\varepsilon)$, 
we obtain
\begin{align}
&\frac1\varepsilon\int_t^{t+\varepsilon}\!\!\int_{\R^3}\zeta_{R_j}\eta_{\ast}(\rho,\mathbf M)\,\dd\boldsymbol{x}\dd\tau
\notag\\
&\,\,\leq
\int_{\R^3}\zeta_{R_j}\eta_{\ast}(\rho_0,\mathbf M_0)\,\dd\boldsymbol{x}
+\int_0^{t+\varepsilon}\!\!\int_{\R^3}\vartheta_\varepsilon
q_{\ast}\cdot\nabla\zeta_{R_j}
\,\dd\boldsymbol{x}\dd\tau
-\int_0^{t+\varepsilon}\!\!\int_{\R^3}\vartheta_\varepsilon\zeta_{R_j}\mathbf M\cdot F\,\dd\boldsymbol{x}\dd\tau.
\notag
\end{align}

For fixed $j$, the Lebesgue differentiation theorem gives
\[
\frac1\varepsilon\int_t^{t+\varepsilon}I_j(\tau)\,\dd\tau\longrightarrow I_j(t)\qquad\text{as }\varepsilon\to 0.
\]
The remaining terms converge by dominated convergence theorem, using the local integrability
of $q_{\ast}$ and \eqref{source-L1}.
We therefore obtain
\begin{align}
\int_{\R^3}\zeta_{R_j}\eta_{\ast}(t)\,\dd\boldsymbol{x}
\leq\int_{\R^3}\zeta_{R_j}\eta_{\ast}(0)\,\dd\boldsymbol{x}
+\int_0^t\!\!\int_{\R^3}q_{\ast}\cdot\nabla\zeta_{R_j}\,\dd\boldsymbol{x}\dd\tau
-\int_0^t\!\!\int_{\R^3}\zeta_{R_j}\mathbf M\cdot F\,\dd\boldsymbol{x}\dd\tau.
\label{fixed-R-relative}
\end{align}

We now let $j\to\infty$. Recall that
$q(\rho,\mathbf M)=\mathbf M\big(\frac{|\mathbf M|^2}{2\rho^2}+h'(\rho)\big)$.
Since $\mathbf M=\rho\mathbf u$, Young's inequality gives
\[
|q(\rho,\mathbf M)|\leq\frac12\rho|\mathbf u|^3+\rho|\mathbf u|h'(\rho)
\leq C\big(\rho|\mathbf u|^3+\rho(h'(\rho))^{3/2}\big).
\]
Moreover, boundedness of $\bar G$ yields
\[
|\bar G\,\mathbf M|\leq C\rho|\mathbf u|\leq C\bigl(\rho+\rho|\mathbf u|^3\bigr).
\]
Consequently,
\[
|q_{\ast}|\leq C\big(\rho|\mathbf u|^3+\rho(h'(\rho))^{\frac{3}{2}}+\rho\big).
\]

By \eqref{flux-condition} and the mass conservation,
$q_{\ast}\in L^1((0,T)\times\R^3)$.
Since
\[
|\nabla\zeta_{R_j}|\leq\frac C{R_j}\mathbf1_{\{R_j<|\boldsymbol{x}|<2R_j\}},
\]
we have
\[
\left|\int_0^t\!\!\int_{\R^3}q_{\ast}\cdot\nabla\zeta_{R_j}\,\dd\boldsymbol{x}\dd\tau\right|
\leq\frac C{R_j}\int_0^T\!\!\int_{R_j<|\boldsymbol{x}|<2R_j}|q_{\ast}|\,\dd\boldsymbol{x}\dd\tau\longrightarrow0.
\]

Furthermore, since $\eta_{\ast}(t)\in L^1(\R^3)$ for almost every $t$, dominated
convergence theorem gives
\[
\int_{\R^3}\zeta_{R_j}\eta_{\ast}(t)\,\dd\boldsymbol{x}\longrightarrow\int_{\R^3}\eta_{\ast}(t)\,\dd\boldsymbol{x}.
\]
The same argument applies to the initial entropy density. Finally, \eqref{source-L1} yields
\[
\int_0^t\!\!\int_{\R^3}\zeta_{R_j}\mathbf M\cdot F\,\dd\boldsymbol{x}\dd\tau
\longrightarrow\int_0^t\!\!\int_{\R^3}\mathbf M\cdot F\,\dd\boldsymbol{x}\dd\tau.
\]

Passing to the limit in \eqref{fixed-R-relative}, using \eqref{int-effective}, and canceling
the constant $\mathcal K_L(\bar\rho)$, we obtain
\[
\mathcal E_L(t)\leq\mathcal E_L(0)-\int_0^t\!\!\int_{\R^3}\mathbf M\cdot F
\,\dd\boldsymbol{x}\dd\tau \qquad\mbox{for {\it a.e.} $t\in(0,T)$}.
\]
Applying \eqref{source-bound} now gives
\[
\mathcal E_L(t)\leq\mathcal E_L(0)+
\begin{cases}
C\displaystyle\int_0^t\bigl(\mathcal E_L(\tau)+\mathcal E_L(\tau)^\chi\bigr)
\,\dd\tau
&\,\,\,\text{in Case 1},\\[3mm]
C_K\displaystyle\int_0^t\mathcal E_L(\tau)\,\dd\tau
&\,\,\,\text{in Case 2}.
\end{cases}
\]

\medskip
\noindent
\textbf{5}. {\it Gr\"onwall estimate and finite-time stability}:  In Case 1, the uniform bound
\eqref{energy-bound} and $\chi\geq1$
imply
\[
\mathcal E_L(\tau)^\chi\leq B_0^{\chi-1}\mathcal E_L(\tau) \qquad\text{for {\it a.e.} }\tau\in(0,T).
\]
Thus the nonlinear source estimate reduces to a linear one.

In Case 2, the source estimate is
already linear in $\mathcal E_L$.

Consequently, in both cases,
\begin{equation}\label{Gronwall-linear}
\mathcal E_L(t)\leq\mathcal E_L(0)+
\begin{cases}
C(1+B_0^{\chi-1})
\displaystyle\int_0^t\mathcal E_L(\tau)\,\dd\tau
&\,\,\,\text{in Case 1},\\[3mm]
C_K\displaystyle\int_0^t\mathcal E_L(\tau)\,\dd\tau
&\,\,\,\text{in Case 2}.
\end{cases}
\end{equation}

Since $\mathcal E_L$ is nonnegative and belongs to $L^\infty(0,T)$,
Gr\"onwall's inequality applies and yields
\[
\mathcal E_L(t)\leq C_T\mathcal E_L(0)\qquad\text{for {\it a.e.} }t\in(0,T),
\]
where $C_T$ has the dependence stated in \eqref{CT-cases}.
This proves \eqref{main-stability}.

\medskip
\noindent
\textbf{6}. {\it Control of the weighted distance}:
Finally, the coercivity estimate\eqref{d1-coercivity} gives
\[
\begin{aligned}
&\int_{\Gamma\cap\{|f|\leq1\}}|f|^2\,\dd\boldsymbol{x}+\int_{\Gamma\cap\{|f|>1\}}
|f|^{\gamma_2}\,\dd\boldsymbol{x}
+\int_{\Gamma^c\cap\{\rho\leq1\}}\rho^{\gamma_1}\,\dd\boldsymbol{x}
+\int_{\Gamma^c\cap\{\rho>1\}}\rho^{\gamma_2}\,\dd\boldsymbol{x}
+\int_{\Gamma^c}\bar{V}\rho\,\dd\boldsymbol{x}\\[1.5mm]
&\,\, 
\leq Cd_1(\rho,\bar\rho).
\end{aligned}
\]
Together with the nonnegativity of $d_2$ and the definition \eqref{EL-def},
this implies that $\mathcal D(t)\leq C\mathcal E_L(t)$.
Combining this with \eqref{main-stability} establishes
\[
\mathcal D(t)\leq C\mathcal E_L(t)\leq CC_T\mathcal E_L(0)
\]
for almost every $t\in(0,T)$, which proves \eqref{distance-bound} and completes the proof.
\end{proof}

\begin{remark}
An early application of the relative entropy method to $L^2$-type stability for multidimensional compressible flows was developed by {\rm Chen--Chen \cite{ChenChen2007}}.
They established stability results for rarefaction waves and vacuum states of the multidimensional Euler equations within a broad class of bounded entropy solutions. For earlier developments of the relative entropy method, we refer to {\rm Dafermos \cite{Dafermos-book}} and the references therein.
\end{remark}

\section{Relative Total Energy}\label{relative-total-energy}

We now relate the positive rotational relative entropy $\mathcal E_L$ to the difference of the total energies of a perturbed state and the rotating equilibrium.

Define the relative total energy by
\[
\mathcal H(\rho,\mathbf M\mid\bar\rho,\bar{\mathbf M}):=\mathcal H(\rho,\mathbf M)
-\mathcal H(\bar\rho,\bar{\mathbf M}).
\]
Since the equilibrium momentum is given by
\[
\bar{\mathbf M}(\boldsymbol{x})=\bar\rho(\boldsymbol{x})\frac{J(m_{\bar\rho}(r(\boldsymbol{x})))}{r(\boldsymbol{x})}\boldsymbol e_\theta,
\]
the angular momentum identity \eqref{swirl-K} yields
\[
\mathcal H(\bar\rho,\bar{\mathbf M})={}\mathcal K_L(\bar\rho)+\int_{\R^3}h(\bar\rho)\,\dd\boldsymbol{x}
+\frac12\int_{\R^3}\bar\rho\bar{\mathcal W}_\alpha\,\dd\boldsymbol{x}.
\]

The following proposition identifies the relative total energy as the sum of the positive modulated entropy and a nonpositive quadratic interaction term.

\begin{proposition}[Relative energy identity]
\label{energy-identity}
Let $(\rho,\mathbf M)$ satisfy the angular momentum compatibility condition
\eqref{compatibility} and have the same
total mass as $\bar\rho$. Assume that both total energies are finite and,
with $f:=\rho-\bar\rho$,
\begin{equation}\label{Q-space}
f\in L^{\frac{6}{6-\alpha}}(\R^3).
\end{equation}
Suppose also that $\mathcal E_L$ is well defined by \eqref{d1-def}, \eqref{d2-def}, and
\eqref{EL-def}{\rm ;} in particular,
\[
\mathcal K_L(\rho)<\infty,\qquad \bar{\Psi}|\rho-\bar\rho|\in L^1(\R^3).
\]
Then
\begin{equation}\label{energy-split}
\mathcal H(\rho,\mathbf M\mid\bar\rho,\bar{\mathbf M})=\mathcal E_L+\mathcal Q_\alpha,
\end{equation}
where
\begin{equation}\label{Qalpha}
\mathcal Q_\alpha:=\frac12\int_{\R^3}f(\Phi_\alpha*f)\,\dd\boldsymbol{x}
\leq0.
\end{equation}
\end{proposition}

\begin{proof}
Since $\rho$ and $\bar\rho$ have the same total mass,
\[
\int_{\R^3}\mu(\rho-\bar\rho)\,\dd\boldsymbol{x}=0.
\]
Set $f:=\rho-\bar\rho$. Using $\bar{\mathcal W}_\alpha
=\Phi_\alpha*\bar\rho$ and the symmetry of the Riesz kernel, we obtain the quadratic expansion
\[
\begin{aligned}
\frac12\int_{\R^3}\rho(\Phi_\alpha*\rho)\,\dd\boldsymbol{x}
-\frac12\int_{\R^3}\bar\rho(\Phi_\alpha*\bar\rho)\,\dd\boldsymbol{x}
=\int_{\R^3}\bar{\mathcal W}_\alpha f\,\dd\boldsymbol{x}
+\frac12\int_{\R^3}f(\Phi_\alpha*f)\,\dd\boldsymbol{x}.
\end{aligned}
\]
All terms are finite under the stated assumptions, by using the
Hardy--Littlewood--Sobolev inequality and the boundedness of
$\bar{\mathcal W}_\alpha$.

Since $\bar G=\bar{\mathcal W}_\alpha+\mu$ and the mass difference vanishes,
\[
\int_{\R^3}\bar G f\,\dd\boldsymbol{x}=\int_{\R^3}\bar{\mathcal W}_\alpha f\,\dd\boldsymbol{x}.
\]
Combining the preceding identities with the rotational cancellation formula
\eqref{EL-effective} proves \eqref{energy-split}.

It remains to establish the sign of $\mathcal Q_\alpha$.
Set
\[
p_\alpha:=\frac6{6-\alpha}.
\]
Then $\frac2{p_\alpha}+\frac{\alpha}{3}=2$.
The Hardy--Littlewood--Sobolev inequality therefore gives
\begin{equation}\label{HLS-quadratic}
\Big|\iint_{\R^3\times\R^3}\frac{g(\boldsymbol{x})h(\boldsymbol{y})}
{|\boldsymbol{x}-\boldsymbol{y}|^\alpha}\,\dd\boldsymbol{x}\dd\boldsymbol{y}\Big|
\leq C_\alpha \|g\|_{L^{p_\alpha}}\|h\|_{L^{p_\alpha}}.
\end{equation}

The Fourier transform of the Riesz kernel
satisfies
\[
\widehat{|\boldsymbol{x}|^{-\alpha}}
(\boldsymbol{\xi})
=c_\alpha|\boldsymbol{\xi}|^{\alpha-3}\qquad \mbox{with $c_\alpha>0$}.
\]
For a real-valued $g\in C_c^\infty(\R^3)$, by Plancherel's identity and the convolution theorem, we obtain
\[
\iint_{\R^3\times\R^3}\frac{g(\boldsymbol{x})g(\boldsymbol{y})}
{|\boldsymbol{x}-\boldsymbol{y}|^\alpha}\,\dd\boldsymbol{x}\dd\boldsymbol{y}
=
\frac{c_\alpha}{(2\pi)^3}\int_{\R^3}|\boldsymbol{\xi}|^{\alpha-3}
|\widehat g(\boldsymbol{\xi})|^2\,\dd\boldsymbol{\xi}\geq0,
\]
where the factor $(2\pi)^{-3}$ corresponds to the Fourier convention.

Now let $f_n\in C_c^\infty(\R^3)$ be real-valued and satisfy $f_n\longrightarrow f$ in $L^{p_\alpha}(\R^3)$.
By \eqref{HLS-quadratic}, the associated quadratic forms converge. Then
\[
\iint_{\R^3\times\R^3}\frac{f(\boldsymbol{x})f(\boldsymbol{y})}{|\boldsymbol{x}-\boldsymbol{y}|^\alpha}\,\dd\boldsymbol{x}\dd\boldsymbol{y}\geq0.
\]

Finally, recalling that $\Phi_\alpha(\boldsymbol{x})
=-\frac1\alpha|\boldsymbol{x}|^{-\alpha},
$
we conclude that
\[
\begin{aligned}
\mathcal Q_\alpha=-\frac1{2\alpha}\iint_{\R^3\times\R^3}\frac{f(\boldsymbol{x})f(\boldsymbol{y})}
{|\boldsymbol{x}-\boldsymbol{y}|^\alpha}\,\dd\boldsymbol{x}\dd\boldsymbol{y}\leq0.
\end{aligned}
\]
This proves \eqref{Qalpha} and completes the proof.
\end{proof}

We next verify that the integrability condition \eqref{Q-space} is satisfied
in both regimes of Theorem \ref{main}.

In Case 1, Lemma \ref{density-norms} applies because
\[
\frac6{6-\alpha}<\frac9{7-2\alpha}\leq\underline\gamma,\qquad 0<\alpha<2.
\]
Consequently,
$
f\in L^{\frac6{6-\alpha}}(\R^3).
$

In Case 2, the prescribed mass and high-integrability assumptions give
\[
f\in L^1(\R^3)\cap L^{\bar p}(\R^3).
\]
Since
\[
1<\frac6{6-\alpha}<\frac3{2-\alpha}<\bar p,
\]
interpolation yields $f\in L^{\frac6{6-\alpha}}(\R^3)$.
Therefore, Proposition \ref{energy-identity} applies in both cases of the main theorem.

\begin{corollary}[Relative total energy inequality]
\label{monotonicity}
Under the assumptions of {\rm Theorem \ref{main}}, suppose additionally in {\rm Case 2} that
$\rho_0\in L^{\bar p}(\R^3)$.
Then
\begin{equation}\label{energy-monotone}
\mathcal H(\rho(t),\mathbf M(t)\mid\bar\rho,\bar{\mathbf M})
\leq\mathcal H(\rho_0,\mathbf M_0\mid\bar\rho,\bar{\mathbf M}) \qquad
\mbox{for \it a.e. $t\in(0,T)$}.
\end{equation}

In particular, both relative energies admit the decomposition \eqref{energy-split}.
\end{corollary}

\begin{proof}
Since $(\bar\rho,\bar{\mathbf M})$ is a stationary rotating equilibrium,
its total energy $\mathcal H(\bar\rho,\bar{\mathbf M})$ is independent of time.

Subtracting this constant from the total energy inequality \eqref{energy-ineq} gives
\[
\mathcal H(\rho(t),\mathbf M(t))-\mathcal H(\bar\rho,\bar{\mathbf M})
\leq\mathcal H(\rho_0,\mathbf M_0)-\mathcal H(\bar\rho,\bar{\mathbf M})
\qquad\mbox{{for \it a.e.}$\, t\in(0,T)$}.
\]
This proves \eqref{energy-monotone}.

The integrability condition \eqref{Q-space} follows for the time-dependent density from the
preceding argument. For the initial density, it follows from
Lemma \ref{density-norms} in Case 1 and from interpolation in Case 2,
using the additional assumption $\rho_0\in L^{\bar p}(\R^3)$.

Proposition \ref{energy-identity} therefore identifies the relative
total energies with $\mathcal E_L+\mathcal Q_\alpha$,
as claimed.
\end{proof}

\begin{remark}
The nonpositivity of $\mathcal Q_\alpha$ follows from the positive definiteness
of the Riesz kernel, rather than merely from the pointwise negativity of
$\Phi_\alpha$. This distinction is essential because the density perturbation
$f=\rho-\bar\rho$ generally changes sign.

Consequently, the relative total energy inequality alone does not provide coercive
control of the perturbation. Indeed, although $\mathcal E_L\geq0$, the
quadratic interaction term $\mathcal Q_\alpha\leq0$ may reduce
the relative total energy.

The positive rotational relative entropy $\mathcal E_L$, together with the vacuum
weight $\bar{V}\rho$ and the nonnegative rotational remainder $d_2$, provides
the weighted $L^2$-type control established in {\rm Theorem~\ref{main}}.
\end{remark}

\section{Rotating White Dwarfs}\label{white-dwarfs}
We now apply the relative entropy stability theory to rotating white dwarfs governed by the
Chandrasekhar equation of state. We specialize to the Newtonian gravitational interaction
$\alpha=1$ and consider the pressure law
\begin{equation}\label{WD-pressure}
p_{\rm wd}(\rho)
=A\int_0^{B\rho^{\frac{1}{3}}}\frac{s^4}{\sqrt{D+s^2}}\,\dd s\qquad \mbox{with $A,B,D>0$}.
\end{equation}
Differentiating gives
\begin{equation}\label{WD-pprime}
p_{\rm wd}'(\rho)=\frac{AB^5}{3}\frac{\rho^{\frac{2}{3}}}{\sqrt{D+B^2\rho^{\frac{2}{3}}}}>0
\qquad \mbox{for $\rho>0$}.
\end{equation}

Let $h_{\rm wd}$ denote the corresponding internal energy density defined by \eqref{h}.
The pressure exhibits different power-law behaviours at low and high densities.
Indeed, direct asymptotic analysis of
\eqref{WD-pressure} yields
\begin{align}
&p_{\rm wd}(\rho)\sim\frac{AB^5}{5\sqrt D}\rho^{\frac{5}{3}},\quad\,\,
h_{\rm wd}(\rho)\sim\frac{3AB^5}{10\sqrt D}\rho^{\frac{5}{3}}
\qquad\,\,\,\mbox{for $\rho\to 0$},
\label{WD-small}\\
&p_{\rm wd}(\rho)
\sim\frac{AB^4}{4}\rho^{\frac{4}{3}},
\quad\,\,
h_{\rm wd}(\rho)
\sim\frac{3AB^4}{4}\rho^{\frac{4}{3}}\qquad\quad\mbox{for $\rho\to\infty$}.
\label{WD-large}
\end{align}
Together with \eqref{WD-pprime}, these asymptotics verify {\rm Assumption \ref{pressure}} with
$\gamma_1=\frac53$, $\gamma_2=\frac43$, and $\underline\gamma=\frac43$.
Moreover,
\[
h_{\rm wd}''(\rho)=\frac{p_{\rm wd}'(\rho)}{\rho}
=\frac{AB^5}{3}\frac{\rho^{-\frac{1}{3}}}{\sqrt{D+B^2\rho^{\frac{2}{3}}}}.
\]
In particular, for every $R>0$,
\[
\inf_{0<\rho\leq R}h_{\rm wd}''(\rho)>0,
\]
consistently with the uniform convexity property \eqref{h-convexity}.

For $\alpha=1$, the mass-critical exponent associated with the attractive interaction is
\[
q_c=\frac{3+\alpha}{3}=\frac43,
\]
which coincides with the high-density adiabatic exponent of the Chandrasekhar
pressure law.

The critical Hardy--Littlewood--Sobolev inequality takes the form
\[
\iint_{\R^3\times\R^3}\frac{\varrho(\boldsymbol{x})\varrho(\boldsymbol{y})}{|\boldsymbol{x}-\boldsymbol{y}|}\,\dd\boldsymbol{x}\dd\boldsymbol{y}\\
\leq C_1^*M^{\frac{2}{3}}\int_{\R^3}\varrho^{\frac{4}{3}}\,\dd\boldsymbol{x},
\]
for every nonnegative density $\varrho$ of mass $M$ for which the
right-hand side is finite.

By \eqref{WD-large},
\[
K_\infty:=\liminf_{\rho\to\infty}\frac{h_{\rm wd}(\rho)}{\rho^{\frac{4}{3}}}
=\frac{3AB^4}{4}.
\]
Consequently, Proposition \ref{coercivity} shows that the fixed-mass coercivity
condition \eqref{EM} holds whenever
\[
\frac{C_1^*}{2}M^{\frac{2}{3}}<\frac{3AB^4}{4},
\]
or equivalently,
\[
M<M_{\rm HLS},\qquad M_{\rm HLS}:=\Big(\frac{3AB^4}{2C_1^*}\Big)^{\frac{3}{2}}.
\]
We remark that $C_1^*$ in $M_{\rm HLS}$ is the best interpolation constant in \cite[Theorem 4]{Carrillo_2018}. 
This HLS critical-mass threshold is a {\it sufficient} condition for coercivity of
the total energy. It is not, however, required for the relative entropy stability
estimate established below.

\begin{corollary}[Relative entropy stability of rotating white dwarfs]
\label{white-dwarf}
Let $\alpha=1$ and let $p=p_{\rm wd}$ be the Chandrasekhar pressure law
\eqref{WD-pressure}.
Let $(\bar\rho,\bar{\mathbf M})$ be a rotating white dwarf equilibrium
of mass $M$ in the sense of {\rm Definition \ref{rotating-star}},
with $L=J^2$ satisfying
\eqref{L-condition}.
Let $(\rho,\mathbf M)$ be an axisymmetric finite-energy entropy weak solution
in the sense of {\rm Definition \ref{entropy-solution}}.
Assume that the solution and its initial data satisfy conditions {\rm(i)} and
{\rm(iii)} of {\rm Theorem \ref{main}}.
Suppose further that, for some $\bar p>3$,
\begin{equation}\label{WD-Lp}
\operatorname*{ess\,sup}_{0<t<T}\|\rho(t)\|_{L^{\bar p}(\R^3)}<\infty.
\end{equation}
Then there exists a constant $C_T>0$, depending on the fixed parameters,
the equilibrium, $T$, and the bound in \eqref{WD-Lp}, such that
\begin{equation}\label{WD-stability}
\mathcal E_L(t) \leq C_T\mathcal E_L(0)
\qquad\text{for {\it a.e.} }t\in(0,T).
\end{equation}
In particular, no restriction $M<M_{\rm HLS}$ is imposed.
\end{corollary}

\begin{proof}
As shown above, the Chandrasekhar pressure law satisfies {\rm Assumption \ref{pressure}} with
$\underline\gamma=\frac{4}{3}$.

For $\alpha=1$, the threshold in Case 2 of Proposition~\ref{weighted-Riesz} becomes
\[
\max\big\{1,\frac6{7-2\alpha}\big\}=\frac65<\frac43=\underline\gamma,
\]
while $\frac3{2-\alpha}=3$.
Thus, the pressure exponent satisfies the required threshold, and
\eqref{WD-Lp} provides the additional high-integrability assumption
with $\bar p>3$.

All hypotheses of Case 2 of {\rm Theorem \ref{main}} are therefore
satisfied. Then applying  {\rm Theorem \ref{main}} yields
\eqref{WD-stability}.
\end{proof}

\bigskip
\medskip
\noindent{\bf Acknowledgments.}
Gui-Qiang G. Chen was partially supported by EPSRC grants  EP/V051121/1 and EP/V008854. 
Difan Yuan was partially supported by EPSRC grants  EP/V051121/1 and
the Fundamental Research Funds for the Central Universities No. 2233100021 and No. 2233300008.
For the purpose of open access, the authors have applied a CC BY public copyright license 
to any Author Accepted Manuscript (AAM) version
arising from this submission.

\bigskip
\noindent{\bf Conflict of Interest:} The authors declare that they have no conflict of interest.
The authors also declare that this manuscript has not been previously published,
and will not be submitted elsewhere before your decision.

\bigskip
\noindent{\bf Data availability:} Data sharing is not applicable to this article as no datasets were generated 
or analyzed during the current study.

\end{document}